\documentclass[aos,MSNbibl,seceqn,nameyear,dvips]{arximspdf}
\usepackage{graphicx}

\doi{10.1214/11-AOS958}
\volume{40}
\issue{1}
\pubyear{2012}
\firstpage{322}
\lastpage{352}

\makeatletter

\renewcommand{\hat}{\widehat}

\renewcommand{\a}{\alpha}
\newcommand{\be}{\beta}
\newcommand{\bX}{{\bar X}}
\newcommand{\bY}{{\bar Y}}
\newcommand{\bigmi}{\mid}
\newcommand{\bigtrip}{|\!|\!|}
\newcommand{\cC}{{\mathcal C}}
\newcommand{\cD}{{\mathcal D}}
\newcommand{\cI}{{\mathcal I}}
\newcommand{\cX}{{\mathcal X}}
\newcommand{\cov}{\operatorname{cov}}
\newcommand{\De}{\Delta}
\newcommand{\ep}{\varepsilon}
\newcommand{\ga}{\gamma}
\newcommand{\hal}{{\hat\a}}
\newcommand{\half}{^{1/2}}
\newcommand{\hbe}{{\hat\be}}
\newcommand{\hg}{{\hat g}}
\newcommand{\hga}{{\hat\ga}}
\newcommand{\hH}{{\widehat H}}
\newcommand{\hK}{{\widehat K}}
\newcommand{\hkb}{{\widehat{K}(b)}}
\newcommand{\hkjb}{{\widehat{K}^j(b)}}
\newcommand{\hkjob}{{\widehat{K}^{j+1}}(b)}
\newcommand{\hkkb}{{\widehat{K}^k(b)}}
\newcommand{\hphi}{{\hat\phi}}
\newcommand{\hpsi}{{\hat\psi}}
\newcommand{\hth}{{\hat\th}}
\newcommand{\inti}{\int_\cI}
\newcommand{\intii}{\inti\inti}
\newcommand{\la}{\lambda}
\newcommand{\mhf}{^{-1/2}}
\newcommand{\mo}{^{-1}}
\newcommand{\mt}{^{-2}}
\newcommand{\mth}{^{-3}}
\newcommand{\mi}{\mid}
\newcommand{\oon}{{1\over n}}
\newcommand{\PC}{^{\mathrm{PC}}}
\newcommand{\pred}{_{\mathrm{pred}}}
\newcommand{\ra}{\to}
\newcommand{\rai}{\ra\infty}
\newcommand{\cent}{^{\mathrm{cent}}}
\newcommand{\si}{\sigma}
\newcommand{\sumion}{\sum_{i=1}^n}
\newcommand{\sumj}{\sum_j}
\newcommand{\sumjoi}{\sum_{j=1}^\infty}
\newcommand{\sumjop}{\sum_{j=1}^p}
\newcommand{\sumk}{\sum_k}
\newcommand{\sumkoi}{\sum_{k=1}^\infty}
\newcommand{\sumkop}{\sum_{k=1}^p}
\newcommand{\sumroi}{\sum_{r=1}^\infty}
\newcommand{\T}{^{\mathrm{T}}}
\newcommand{\tg}{{\tilde g}}
\renewcommand{\th}{\theta}
\newcommand{\tiH}{{\widetilde H}}
\newcommand{\trip}{|\!|\!|}
\newcommand{\var}{\operatorname{var}}
\newcommand{\ze}{\zeta}

\newtheorem{theo}{Theorem}[section]

\newproclaim{rem}{Remark}



\begin{document}
\begin{frontmatter}

\title{Methodology and theory for partial least squares applied to functional data\thanksref{T1}}
\runtitle{Partial least squares}

\thankstext{T1}{Supported by grants and fellowships
from the Australian Research Council.}

\begin{aug}
\author[A]{\fnms{Aurore} \snm{Delaigle}\corref{}\ead[label=e1]{A.Delaigle@ms.unimelb.edu.au}}
\and
\author[B]{\fnms{Peter} \snm{Hall}\ead[label=e2]{halpstat@ms.unimelb.edu.au}}
\runauthor{A. Delaigle and P. Hall}
\affiliation{University of Melbourne, and University of Melbourne
and~University~of~California, Davis}
\address[A]{Department of Mathematics and Statistics\\
University of Melbourne\\
Parkville, VIC, 3010\\
Australia\\
\printead{e1}}
\address[B]{Department of Statistics\\
University of California\\
Davis, California 95616\\
USA\\
\printead*{e2}} %adresu isvedimo komanda gale!
\end{aug}

% HISTORY:
\received{\smonth{8} \syear{2011}}
\revised{\smonth{11} \syear{2011}}

% ABSTRACT
%
\begin{abstract}
The partial least squares procedure was originally developed to
estimate the slope parameter in multivariate parametric models. More
recently it has gained popularity in the functional data literature.
There, the partial least squares estimator of slope is either used to
construct linear predictive models, or as a~tool to project the data
onto a~one-dimensional quantity that is employed for further
statistical analysis. Although the partial least squares approach is
often viewed as an attractive alternative to projections onto the
principal component basis, its properties are less well known than
those of the latter, mainly because of its iterative nature. We
develop an explicit formulation of partial least squares for functional
data, which leads to insightful results and motivates new theory,
demonstrating consistency and establishing convergence rates.
\end{abstract}

% KEYWORDS
%

\end{frontmatter}

%s1 #&#
\section{Introduction}

Partial least squares (PLS) is an iterative procedure for estimating
the slope of linear models. The technique was originally developed in
high-dimensional and collinear multivariate settings and is especially
popular in chemometrics. See \citet{Wol75}, \citet{MarNs89},
\citet{Hel90},
\citet{FraFri93}, \citet{Gar94}, \citet{GouFea96},
\citet{DurSab97}
and \citet{NguRoc04}.

The iterative nature of PLS can make it difficult to uncover properties
in a~clear and explicit way, and for a~long time PLS was regarded as a~technique that worked well, but whose properties were relatively
obscure. Early theoretical developments of multivariate PLS can be
found in \citet{LorWanKow87} and \citet{Hos88},
and further developments include those of \citet{PhaReiPen02},
\citet{PhadeH03}, \citet{BroEld09} and
\citet{KraSug11}.

More recently, the method has been applied in the functional data
context by \citet{PreSap05N1}, who suggest using PLS for
estimating slope in functional linear models; see also \citet
{ReiOgd07}. Also in the functional setting, the intrinsic iterative nature
of PLS has made it difficult to develop intuition and derive clear and
explicit theoretical properties. In this paper we provide a~transparent
account of theoretical issues that underpin PLS methods in linear
models for prediction from functional data, and show that they motivate
an alternative formulation of PLS in that setting. This ``alternative
PLS,'' which we refer to here as APLS, has the advantage that it is
expressed only in terms of functions that are explicitly computable.
These attributes make APLS particularly attractive, relative to the
conventional PLS formulation, and permit detailed theoretical development.

We give concise stochastic expansions for the difference between
estimators derived using APLS, and the quantities to which these
estimators converge in the large-sample limit. These expansions are
valid uniformly in estimators based on the first $O(n\half)$ APLS
basis functions, where $n$ denotes sample size. The expansions also
lead easily and directly to a~variety of results about our estimators,
including convergence rates and central limit theorems.

Besides functional linear models, PLS is employed in a~variety of other
data functional problems.
For example, \citet{FerVie06} use it to define a~semi-metric for
nonparametric functional predictors or classifiers;
\citet{EscAguVal07} employ PLS with logit regression;
\citet{PreSapLev07} and \citet{DelHal12} use it for
functional data classification. See also \citet{PreSap05N2},
\citet{KraBouTut08} and \citet{Aguetal10}.

%s2 #&#
\section{Functional linear models}
%s2.1 #&#
\subsection{General bases for inference in functional linear
models}\label{sec22}

Let $\cX=\{(X_1,Y_1),\ldots,(X_n,Y_n)\}$ denote a~sample of
independent data pairs, all distributed as $(X,Y)$, where $X$ is a~random function defined on the nondegenerate, compact interval $\cI$
and satisfying $\inti E(X^2)<\infty$, and $Y$ is a~scalar random
variable generated by the linear model
%
%e2.1 #&#
%
\begin{equation}\label{21}
Y=a+\inti b X+\ep.
\end{equation}
Here, $a$ denotes a~scalar parameter, $\ep$ is a~scalar random
variable with finite mean square and satisfying $E(\ep\mi X)=0$ and
$b$, a~function-valued parameter, is a~square-integrable function
on $\cI$.

Predicting the value of $Y$, given $X$, amounts to estimating the function
%
%e2.2 #&#
%
\begin{equation}\label{22}
g(x)=E(Y\mi X=x)=a+\inti b x ,
\end{equation}
which, itself, requires us to estimate the scalar $a$ and the function
$b$ from the data. A standard approach is to express $X$ and $b$ in
terms of an orthonormal basis $\psi_1,\psi_2,\ldots$ defined on $\cI
$. Expansions for $X$ and $b$ in this basis can be written as $X=\sumj
(\inti X \psi_j) \psi_j$ and $b=\sumj v_j \psi_j$, where
$v_j=\inti b \psi_j$.
Since, in practice, we can calculate only a~finite number of terms, the
infinite-dimensional expansion for $b$ is approximated by a~sum of $p$
terms, where $p\geq1$ is an integer, and each term of this sum is then
estimated from the data.
Note that $\inti b X=\sumj v_j\inti X \psi_j$, which motivates us to
take $a=E(Y)-\inti b E(X)$ and define $\be_1,\ldots,\be_p$ to be
the sequence $v_1,\ldots,v_p$ that minimizes
%
%e2.3 #&#
%
\begin{equation}\label{23}
s_p(v_1,\ldots,v_p)
=E\Biggl\{\inti b (X-EX)-\sumjop v_j\inti(X-EX) \psi_j\Biggr\}^{
2} .
\end{equation}
The functions
%
%e2.4 #&#
%
\begin{eqnarray}\label{24}
b_p&=&\sumjop\be_j\psi_j,\nonumber\\[-8pt]\\[-8pt]
g_p(x)&=&E(Y)+\inti b_p (x-EX)=E(Y)+\sumjop
\be_j\inti(x-EX)\psi_j\nonumber
\end{eqnarray}
are approximations to $b$ and to $g(x)$, respectively. Their accuracy,
as $p$ increases, depends on the choice of the sequence $\psi_1,\psi
_2,\ldots.$

Sometimes the basis is chosen independently of the data
(e.g., sine-cosine basis, spline basis, etc). Then the functions $\psi
_j$ are known, and an empirical version of (\ref{24}) is obtained
by replacing the scalars $\beta_1,\ldots,\beta_p$ by the sequence
$v_1,\ldots,v_p$ that minimizes
%
%e2.5 #&#
%
\begin{equation}\label{LSconv}
n^{-1} \sumion\Biggl\{Y_i-\bar Y-\sumjop v_j \inti(X_i-\bX) \psi
_j\Biggr\}^{2}.
\end{equation}
A drawback of such bases is that there is no reason why their first $p$
elements should capture the most important information about the
regression function $g$, available from the data. It seems more
attractive to use bases that adapt to the properties of the population
represented by the data. We discuss two such adaptive bases in Sections
\ref{sec23} and~\ref{sec24}, respectively.

%s2.2 #&#
\subsection{Principal component basis}\label{sec23}
One of the most popular adaptive bases is the so-called principal
component basis, constructed from the covariance function\vadjust{\goodbreak} $K(s,t)=\cov
\{X(s),X(t)\}$ of the random process $X$. As is common in
mathematical analysis, we shall use the notation $K$ also for the
linear transformation (a~functional) that takes a~square-integrable
function~$\psi$ to~$K(\psi)$ given by
$
K(\psi)(t)=\inti\psi(s) K(s,t) \,ds .
$

Since $\inti E(X^2)<\infty$, then $\inti K(t,t) \,dt<\infty$, and we
can write the spectral decomposition of $K$ as
%
%e2.6 #&#
%
\begin{equation}\label{25}
K(s,t)=\sumkoi\th_k \phi_k(s) \phi_k(t) ,
\end{equation}
where the principal component basis $\phi_1,\phi_2,\ldots$ is a~complete orthonormal sequence of eigenvectors (i.e., eigenfunctions) of
the transformation $K$, with respective nonnegative eigenvalues $\th
_1,\th_2,\ldots.$ That is, $K(\phi_k)=\th_k \phi_k$ for \mbox{$k\geq
1$}. Positive definiteness of $K$ implies that the eigenvalues are
nonnegative, and the condition $\inti E(X^2)<\infty$ entails $\sumk
\th_k<\infty$. Therefore we can, and do, order the terms in the
series in (\ref{25}) so that
%
%e2.7 #&#
%
\begin{equation}\label{26}
\th_1\geq\th_2\geq\cdots\geq0 .
\end{equation}

In practice the scalars $\th_j$ and the functions $\phi_j$ are
unknown and are estimated from the data, as follows. First, the
covariance function is estimated by
%
%e2.8 #&#
%
\begin{equation}\label{cov}
\hK(s,t)=\oon\sumion\{X_i(s)-\bX(s)\} \{X_i(t)-\bX(t)\},
\end{equation}
where\vspace*{2pt} $\bX(t)=n^{-1}\sumion X_i(t)$. Then, $\th_1,\ldots,\theta_n$
and $\phi_1,\ldots,\phi_n$ are estimated by the eigenvalues $\hth
_1\geq\hth_2\geq\cdots\hth_n\geq0$ and the eigenfunctions $\hphi
_1,\ldots,\hphi_n$ of the transformation represented by $\hat K$,
which can have at most $n$ nonzero eigenvalues.
Finally, an empirical version of $\beta_1,\ldots,\beta_p$ is defined
to be the sequence $v_1,\ldots,v_p$ that minimizes (\ref{LSconv}),
where each~$\psi_j$ there is replaced by~$\hat\phi_j$. Then,~$g_p$
at (\ref{24}) is replaced by its corresponding empirical version.
In the rest of this paper, to avoid confusion with projections of $b$
onto other bases, we shall add a~superscript $\PC$ to coefficients
obtained from projection of $b$ onto one of the functions $\phi_j$;
that is, we shall use the notation $\beta_j\PC=\inti
b\phi_j$.\looseness=-1

The literature on functional linear models based on principal component
analysis (PCA) is large. It includes, for example, work by \citet
{CaiHal06}, \citet{ReiOgd07}, \citet{ApaGol08},
\citet{CarSar08}, \citet{Bal09}, \citet{MulYao10},
\citet{WuFanMul10} and \citet{YaoMul10}.

%s2.3 #&#
\subsection{The orthonormal PLS basis}\label{sec24}
The principal component basis introduced in Section~\ref{sec23} is
defined in terms of the population, but only through~$X$. In
particular, while its first $p$ elements $\phi_1,\ldots,\phi_p$
usually contain most of the information related to the covariance of
$X$, these are not necessarily important for representing $b$, and all
or some of the most important terms accounting for the interaction
between $b$ and $X$ might come from later principal components. In
prediction, to capture the main effects of interaction using only a~few
terms, one could construct the basis in a~way that takes this
interaction into account.

Motivated by such considerations, the standard PLS basis, adapted to
the functional context, is defined iteratively by choosing $\psi_p$ in
a~sequential manner, to maximize the covariance functional
%
%e2.9 #&#
%
\begin{equation}\label{27}
f_p(\psi_p)=\cov\biggl\{Y-g_{p-1}(X),\inti X \psi_p\biggr\},
\end{equation}
subject to
%
%e2.10 #&#
%
\begin{eqnarray}\label{28b}
\intii\psi_j(s)K(s,t) \psi_p(t) \,ds \,dt&=&0 \qquad\mbox{for }
1\leq j\leq p-1\quad\mbox{and}\nonumber\\[-8pt]\\[-8pt]
\|\psi_p\|&=&1 ,\nonumber
\end{eqnarray}
where \mbox{$\|\cdot\|$} is a~norm (see Section~\ref{sec24b}), and given
that $\psi_1,\ldots,\psi_{p-1}$ have already been chosen. [Recall
that $g_p$ was defined at (\ref{24}).]
In practice, the covariances in (\ref{27}) are replaced by
estimates, and empirical versions of the~$\psi_j$'s are constructed by
an iterative algorithm described in Appendix~\ref{sec31}.

Partial least squares can also be used for prediction in nonlinear
models, where the basis that it produces is sometimes, but not always,
effective for prediction. Specifically, although the PLS basis enables
a~consistent approximation to $g$ in such cases, a~large number of
terms may be required to get a~good approximation.

%s3 #&#
\section{Properties of theoretical functional partial least
squares}\label{sec2}

For prediction and estimation of $b$, the PLS basis is sometimes
preferred to the PCA basis, partly because it can often capture the
relevant information with fewer terms; see our data illustrations in
Section~\ref{sec4}. Detailed theoretical properties for inference in
functional linear models based on the PCA basis have been studied in a~number of papers, but few results exist about their functional PLS counterpart.
In this section we provide new insight into the theoretical PLS basis,
defined in (\ref{27}) and (\ref{28b}), and give an explicit
description of the space generated by the first $p$ PLS basis functions
$\psi_1,\ldots,\psi_p$. These properties motivate an alternative
formulation of functional PLS, which we call APLS. It permits us to
define the functional PLS basis very simply, and to construct an
explicitly defined algorithm to implement empirical PLS; see Section
\ref{sec3}. The explicit nature of the algorithm will allow us to
derive detailed theoretical properties of empirical functional PLS,
including convergence rates; see Section~\ref{sec5}.

%s3.1 #&#
\subsection{Explicit form of the orthonormal PLS basis}\label{sec24b}
Our first result, Theorem~\ref{Theorem1}, below, gives an
explicit\vadjust{\goodbreak}
account of the constrained optimization described in Section
\ref{sec24}. We use the following notation. Given $\a_1$ and $\a_2$ in
the class $\cC(\cI)$ of all square-integrable functions\vspace*{1pt} on
$\cI$, write $\intii\a_1 \a_2 K$ to denote $ \intii\a_1(s) \a_2(t)
K(s,t) \,ds \,dt$. For any\vspace*{1pt} $x\in\cC(\cI)$, define
$\|x\|^2=\intii x x K$. (Some implementations of PLS, e.g., the one in
Appendix~\ref{sec31}, take $\|x\|^2=\inti x^2$, but this affects only
the scale, not the main properties of the functions $\psi_j$.)
%
%th3.1 #&#
%
\begin{theo}\label{Theorem1} If $\inti E(X^2)<\infty$, then the
function $\psi_p$ that maximizes~$f$ at (\ref{27}), given $\psi
_1,\ldots,\psi_{p-1}$ and subject to %\eqref{28} or
(\ref{28b}), is determined by
%
%e3.1 #&#
%
\begin{equation}\label{29}
\psi_p=c_0 \Biggl[K\Biggl\{b-\sum_{j=1}^{p-1}
\biggl(\inti b \psi_j\biggr) \psi_j\Biggr\}
+\sum_{k=1}^{p-1} c_k \psi_k\Biggr] ,
\end{equation}
where,
for $1\leq k\leq p-1$, the constants $c_k$ are obtained by solving the
linear system of $p-1$ equations
%
%e3.2 #&#
%
\begin{equation}\label{210b}
\intii\psi_j \psi_p K=0,\qquad j=1,\ldots,p-1,
\end{equation}
and where
$c_0$ is defined uniquely, up to a~sign change, by the property
%
%e3.3 #&#
%
\begin{equation}\label{211}
\|\psi_p\|=1 .
\end{equation}
\end{theo}

One of the interesting implications of the theorem is that for each
$j$, the $j$th basis function determined by PLS can be expressed as a~linear combination of $j$ explicitly defined functions. More precisely,
the theorem implies that $\psi_1=d_1 K(b)$, where, by (\ref{211})
with $p=1$, $d_1=\|K(b)\|\mo$, and more generally, the following
properties follow from the representation (\ref{29}); the first
property implies the second:
%
%e3.4 #&#
%
\begin{equation}
\label{212}
\begin{tabular}{p{320pt}}
(a) For each $p\geq1$, and given $\psi_1,\ldots,\psi_{p-1}$, the
function $\psi_p$ is the linear
combination of $K(b),\ldots,K^p(b)$ for which %\eqref{28} or
(\ref{28b}) holds, and is unique up
to a~sign change. (b) For each $p\geq1$, representing a~function as a~linear
form in $\psi_1,\ldots,\psi_p$ is equivalent to representing it as a~linear combination
of $K(b),\ldots,K^p(b)$.
\end{tabular}\hspace*{-28pt}
\end{equation}
These properties motivate the APLS formulation and underpin the rest of
the paper.
Interestingly, (\ref{212}) continues to hold if equations (\ref
{210b}) are replaced by $\inti\psi_j\psi_p=0$ for $j=1,\ldots
,p-1$. In particular, although the functions $\psi_2,\ldots,\psi_p$
will change in this case, the spaced spanned by $\psi_1,\ldots,\psi
_p$ will not alter.

Result (\ref{212}) is a~deterministic functional version of a~known
result for empirical PLS in the multivariate context. More
specifically, suppose we have $n$ observations of a~$q$-variate
(nonfunctional) predictor of a~variable $Y$, let $\mathbf X$ be the
$n\times q$ matrix of observations on the predictor, and let
${\mathbf y}$ be the $n\times1$ vector containing the
observations on $Y$. Then it has been established that the space
spanned by the first $p$ empirical PLS components is equal to the space
generated by ${\mathbf X}\T{\mathbf y},{\mathbf
D}{\mathbf X}\T{\mathbf y},\ldots,{\mathbf
D}^{p-1}{\mathbf X}\T{\mathbf y}$, where ${\mathbf
D}={\mathbf X}\T{\mathbf X}$. See, for example,
\citet{BroEld09}, and compare the empirical algorithm in Section
\ref{sec32}. This is itself a~particular case of results that are
available more generally in Krylov spaces, although again in the
multivariate rather than functional setting, that is the subject of
this paper.
%Krylov space

%s3.2 #&#
\subsection{Expansions in a~nonorthogonal PLS basis}\label{sec25}
The properties at (\ref{212}) give a~clear and explicit account of
the form taken by the PLS basis functions. For example, they show that
for each $p$, the space generated by $\psi_1,\ldots,\psi_p$ is the
same as the space generated (i.e., spanned) by $K(b),\ldots,K^p(b)$.
Note that the functions $K^j(b)$ are explicitly defined, since we have
$K^j(b)=\sum_k\th_k^j \be_k\PC\phi_k$, where $\phi_k$ is the
$k$th PCA basis function.

Next, if we note that $a=E(Y)-\inti b E(X)$ and define $\ga_1,\ldots
,\ga_p$ to be the sequence $w_1,\ldots,w_p$ that minimizes
%
%e3.5 #&#
%
\begin{equation}\label{213}\qquad
t_p(w_1,\ldots,w_p)
=E\Biggl\{\inti(X-EX) b-\sumjop w_j\inti(X-EX) K^j(b)\Biggr\}^{2}
\end{equation}
[compare (\ref{23})],
then the slope function approximation $b_p$ at (\ref{24}) has two
equivalent expressions,
%
%e3.6 #&#
%
\begin{equation}\label{214}
b_p=\sumjop\ga_j K^j(b)=\sumjop\be_j \psi_j ,
\end{equation}
where $\be_1,\ldots,\be_p$ are as defined in Section~\ref{sec22}
if we take the general $\psi_1,\ldots,\psi_p$ introduced there to be
the specific functions given by Theorem~\ref{Theorem1}.

In matrix notation,
%
%e3.7 #&#
%
\begin{equation}\label{56}
\ga\equiv(\ga_1,\ldots,\ga_p)\T=H\mo(\a_1,\ldots,\a_p)\T,
\end{equation}
where $H=(h_{jk})_{1\leq j,k\leq p}$ denotes a~$p\times p$ matrix,
%
%e3.8 #&#
%e3.9 #&#
%
\begin{eqnarray}
\label{54}
h_{jk}&=&\inti K^{j+1}(b) K^k(b)
=\sumroi(\be_r\PC)^2 \th_r^{j+k+1} ,\\
\label{55}
\a_j&=&\inti K(b) K^j(b)=\sumroi(\be_r\PC)^2 \th_r^{j+1}=h_{0j}.
\end{eqnarray}
Here we have used the fact that, for $p$ fixed, the matrix $H$ is
nonsingular because, for finite $p$, the equivalence of the expansion
in the orthogonal basis $\psi_1,\ldots,\psi_p$ and in the basis
$K(b),\ldots,K^p(b)$ implies that the sequence $\ga_1,\ldots,\ga_p$
that minimizes (\ref{213}) is unique. See also our discussion on
Hankel matrices in Section~\ref{secfurther}.

The $p$th order approximation $g_p(x)$ to $g(x)=E(Y\mi X=x)$, resulting
from the $p$th order approximation of $b$ by either of the identities
at (\ref{214}), is given equivalently by the second formula at
(\ref{24}) or by the expression
%
%e3.10 #&#
%
\begin{equation}\label{215}
g_p(x)=a+\inti b_p x=E(Y)+\sumjop\ga_j \inti(x-EX) K^j(b).
\end{equation}
We denote by APLS the formulation of PLS based on the sequence
$K(b),\ldots,\break K^p(b)$.

For the approximation at (\ref{214}) to converge to $b$, that
function should be expressible as a~linear form in
$K(b),K^2(b),\ldots,$
%
%e3.11 #&#
%
\begin{equation}\label{216}
b=\sumjoi w_j K^j(b) ,
\end{equation}
where the $w_j$'s are constants, and the series converges in $L^2$. The
next theorem gives conditions under which, for a~general $b$ in $\cC
(\cI)$, there exist $w_1,w_2,\ldots$ such that (\ref{216}) holds.
%
%th3.2 #&#
%
\begin{theo}\label{Theorem2} If $\inti E(X^2)<\infty$, and the
eigenvalues of $K$ are all nonze\-ro, then each $b\in\cC(\cI)$ can be
written as at (\ref{216}),
where the series converges in~$L^2$.
\end{theo}

Under the side condition $\inti E(X^2)<\infty$ the assumption in
Theorem~\ref{Theorem2} that all eigenvalues of $K$ be nonzero is both
necessary and sufficient for (\ref{216}) to hold for all $b\in\cC
(\cI)$. However, if some eigenvalues $\th_j$, corresponding to
respective eigenvectors $\phi_j$, vanish, then the respective values
of $\inti(X-EX) \phi_j$ vanish with probability 1, and so those
indices make zero contribution to $\inti(X-EX) b=\sumj\inti(X-EX)
\phi_j\cdot\inti b \phi_j$. Therefore we can delete the components
of $b=\sumj\phi_j\inti b \phi_j$ that correspond to indices $j$ for
which $\th_j=0$, without affecting the value of $\inti b X$; and it
is only through the latter integral that $b$ influences prediction.
Therefore the theorem can be stated in a~form which asserts that even
if some of the eigenvalues of $K$ vanish, the representation at (\ref
{216}) is sufficiently accurate for all purposes of prediction based
on (\ref{21}). The only reason we have not taken this course is to
make our arguments relatively simple and transparent.

Note that the $w_j$'s in (\ref{216}) are not determined uniquely. In
particular, (\ref{216}) implies that $K(b)=\sumj w_j K^{j+1}(b)$,
and so the following expansion, among many others, is an alternative
to (\ref{216}):
$b=\sumjoi(w_j+w_{j+1}) K^{j+1}(b)$. %\label{217}
This lack of uniqueness does not violate the equivalence noted in
(\ref{212})(b), since that property is asserted only for a~finite
sequence $\psi_1,\ldots,\psi_p$.
However, it makes it impossible to treat usefully the relationship
between the infinite expansion of a~function $b$ in terms of the
sequence $K(b),K^2(b),\ldots,$ and its infinite expansion in terms of
the PCA basis, $\phi_1,\phi_2,\ldots,$ introduced in Section~\ref{sec23}. Nevertheless we can discuss the $p$th order PLS projection
$b_p=\sumjop\ga_j K^j(b)$ of $b$ onto the finite-dimensional space
spanned by $K(b),\ldots,K^p(b)$, for an arbitrary but fixed $p\geq1$.

To this end, recall that $\be_1\PC,\be_2\PC,\ldots$ denote the
Fourier coefficients of $b$ with respect to the PCA basis $\phi_1,\phi
_2,\ldots.$ Then,
%
%e3.12 #&#
%
\begin{equation}\label{219}\qquad
b_p=\sumjop\ga_j K^j(b)=\sumjop\ga_j \sumkoi\be_k\PC\th
_k^j \phi_k
=\sumkoi\be_k\PC\Biggl(\sumjop\ga_j \th_k^j\Biggr) \phi_k,
\end{equation}
and the last series expresses $b_p$ in terms of the components of the
PCA basis.

%s4 #&#
\section{Empirical implementation of APLS}\label{sec3}

%s4.1 #&#
\subsection{Algorithm for empirical APLS}\label{sec32}
A standard algorithm for empirical implementation of PLS based on the
sequence $\psi_1,\ldots,\psi_p$ is given in Appendix~\ref{sec31}.
In this section we describe a~simple empirical algorithm for
implementing APLS based on the nonorthogonal sequence $K(b),\ldots
,K^p(b)$. As we shall see, this algorithm will permit simple derivation
of theoretical properties of PLS. In Section~\ref{secaltAPLS} we
shall deduce two algorithms that are numerically more stable.

To estimate $K(b),\ldots,K^p(b)$, first note that we can estimate
$K(b)$ by
\[
\hkb=\oon\sumion X_i\cent Y_i\cent
=\oon\sumion(X_i-\bX) (Y_i-\bY) ,
\]
where $X_i\cent=X_i-\bX$ and $Y_i\cent=Y_i-\bY$. Then, given an
estimator $\hkjb$ of $K^j(b)$, we can estimate\vspace*{1pt}
$K^{j+1}(b)(t)$ by $ \hkjob(t)=\inti\hkjb(s) \hK(s,t) \,ds$,  where
$\hK$ is the conventional estimator of the covariance function, $
\hK(s,t)=n\mo\sumion\{X_i(s)-\bX(s)\} \{X_i(t)-\bX(t)\}$.
Having\vspace*{1pt} calculated $\hkjb$ for $1\leq j\leq p$ we take
$\hga_1,\ldots ,\hga_p$ to minimize
%
%e4.1 #&#
%
\begin{equation}\label{34}
U_p(w_1,\ldots,w_p)
=\oon\sumion\Biggl\{Y_i\cent-\sumjop w_j \inti X_i\cent\hkjb
\Biggr\}^{ 2}
\end{equation}
with respect to $w_1,\ldots,w_p$ [compare (\ref{213})]. In matrix notation,
%
%e4.2 #&#
%
\begin{equation}\label{51}
\hga\equiv(\hga_1,\ldots,\hga_p)\T=\hH\mo(\hal_1,\ldots,\hal
_p)\T,
\end{equation}
where $\hH=({\hat h}_{jk})_{1\leq j,k\leq p}$ denotes a~$p\times p$ matrix,
%
%e4.3 #&#
%e4.4 #&#
%
\begin{eqnarray}\qquad
\label{52}
{\hat h}_{jk}&=&\intii\hK(s,t) \hkjb(s) \hkkb(t) \,ds \,dt
=\inti\hkjob\hkkb,\\
\label{53}
\hal_j&=&\inti\hkb\hkjb.
\end{eqnarray}

Finally we construct an estimator of $g$ based on (\ref{215}),
%
%e4.5 #&#
%
\begin{equation}\label{35}
\hg_p(x)=\bY+\sumjop\hga_j \inti(x-\bX) \hkjb.
\end{equation}

%re1 #&#
%
\begin{rem}
Formula (\ref{54}) demonstrates that the theoretical version $H$ of
$\hH$ is a~symmetric matrix. Our estimator $\hH$ does not necessarily
enjoy that property, but an alternative estimator of $h_{jk}$ can be
defined to satisfy it. More precisely we can take
$
{\tilde h}_{jk}=\inti\widehat{K}^{j+k}(b) \hkb,
$
which produces a~symmetric estimator $\tiH=({\tilde h}_{jk})$ of~$H$.
We could use $\tiH$ in place of $\hH$, but computing ${\tilde h}_{jk}$
requires $\hK$ to be iterated $j+k$ times, whereas ${\hat h}_{jk}$
needs iteration\vspace*{2pt} at most $\max(j+1,k)$ times. Therefore we prefer the
version $\hH$.
\end{rem}

%s4.2 #&#
\subsection{Stabilized algorithm for empirical APLS}\label{secaltAPLS}
The algorithm described in Section~\ref{sec32} would provide a~good
solution if we were able to work in exact arithmetic, but it can be
unstable in finite precision arithmetic. This is because, due to the
nonunicity of the expression for $b$ in terms of the infinite series
$K(b),K^2(b),\ldots,$ as $p$ increases the linear system of equations
given by the empirical version of (\ref{213}) [see (\ref{34})]
becomes closer to singular. Therefore, in finite precision arithmetic,
as $p$ increases it becomes more difficult to numerically identify one
or more of the valid expressions arising from a~large number of terms
in the sequence $\hat K(b), \hat K^2(b),\ldots.$

There exist a~number of numerical methods for overcoming this numerical
difficulty. A simple approach is to transform the linear system of
equations by Gram--Schmidt orthogonalization; see Section 7.7 of
\citet{Lan99}. There,\vspace*{1pt} the columns of the $n\times p$ matrix with $(i,j)$th
element equal to $\inti X_i\cent\hkjb$ are transformed into $p$
orthonormal vectors $u_1,\ldots,u_p$ by the modified Gram--Schimdt
algorithm (a~numerically stabilized version of Gram--Schmidt algorithm;
see Appendix~\ref{secgram}).
Instead of using $\hga$ in (\ref{51}), the sequence that minimizes
(\ref{34}) can then be computed by solving, with respect to
$w_1,\ldots,w_p$, the equivalent equation
$
{\mathbf R} (w_1,\ldots,w_p)\T={\mathbf U}\T(Y_1\cent
,\ldots,Y_n\cent)\T,
$
where ${\mathbf U}$ is a~matrix with columns $u_1,\ldots,u_p$, and
${\mathbf R}$ is an upper $p\times p$ triangular matrix.
Let $\hga^*=(\hga_1^*,\ldots,\hga_p^*)\T$ be the solution of this
equation. We can estimate $g$ by
$\hg_p^*(x)=\bY+\sumjop\hga_j^* \inti(x-\bX) \hkjb$. %

Alternatively, having constructed $\hkjb$ for $1\leq j\leq p$ as in
Section~\ref{sec32}, we can also transform them into an orthonormal
sequence $\hpsi_1,\ldots,\hpsi_p$ satisfying the standard PLS
constraints, $\intii\hpsi_j \hpsi_k \hat K=0$ for $j\neq
k$ [(compare~(\ref{28b})], using, for example, the modified
Gram--Schmidt algorithm.
Then we can calculate an empirical version $\hat\beta_1,\ldots,\hat
\beta_p$ of $\beta_1,\ldots,\beta_p$, the latter defined in
Section~\ref{sec22} (taking there the $\psi_j$'s to be the empirical
PLS basis functions), by finding the sequence $v_1,\ldots,v_p$ that
minimizes (\ref{LSconv}).\vadjust{\goodbreak}
Finally, we can estimate~$g$ by
%
%e4.6 #&#
%
\begin{equation}\label{33APLS}
\tg_p(x)=\bY+\sumjop\hbe_j\inti(x-\bX) \hpsi_j .
\end{equation}

In exact arithmetic, $\hg_p^*$ and $\hga^*$ would be equal to,
respectively, $\hg_p$ and $\hga$ defined in (\ref{35}) and (\ref
{51}). Likewise, $\tg_p$, would be equal to $\hg_p$. In practice,
these approximations differ because we can only work in finite
precision arithmetic, and the algorithms leading to $\hg_p^*$ and $\tg
_p$ are much more numerically stable than the one leading to $\hg_p$.
In general, for prediction we found the algorithm leading to $\tg_p$
to be preferable.
However, the algorithm of Section~\ref{sec32} is important for
developing intuition and assembling theoretical arguments. On the
theoretical side, the simple, explicit formulae in Section~\ref{sec32} permit us to establish consistency and derive rates of
convergence. Of course, the equivalence between $\tg_p$, $\hg_p$ and
$\hg_p^*$ implies that, in order to derive the theoretical properties
of $\tg_p$ and $\hg_p^*$, it suffices\vspace*{1pt} to derive them for~$\hg_p$
(all three have the same theoretical properties). On the intuitive side
we note that the explicit formulation of the quantities involved in our
empirical algorithms for APLS gives a~much clearer account of what
partial least-squares does, than the standard empirical iterative PLS
algorithm in Appendix~\ref{sec31}.

%s5 #&#
\section{Asymptotic properties of empirical APLS}\label{sec5}

%s5.1 #&#
\subsection{Introduction}\label{sec51}
To our knowledge, the only existing theoretical results for functional
PLS are those of \citet{PreSap05N1}, who state
generalizations
to the functional data context of some results of \citet{Hos88}.
Although they are of interest, the theoretical arguments there are
iterative and not explicit, and consistency of the PLS approximation is
mentioned without a~proof and without regularity conditions or
convergence rates. This is because those results are based on the
iterative empirical approximation of PLS, and the inexplicit form of
the algorithm (see Appendix~\ref{sec31}) apparently makes it very
difficult to derive explicit theoretical results.

Our alternative formulation, APLS, of the functional partial
least-squares problem permits us to derive many properties. As already
explained in Section~\ref{secaltAPLS}, the theoretical properties of
the empirical approximations~$\hg_p^*$ and~$\tg_p$ in Section~\ref{secaltAPLS} are identical to those of $\hg_p$ in Section~\ref{sec32}.

%s5.2 #&#
\subsection{Main results}\label{sec52}

Define $\mu=E(X)$, a~function, and observe that we can write:
%
%e5.1 #&#
%
\begin{equation}\label{58}\qquad
\hK=K+n\mhf\xi+n\mo\eta,\qquad \hkb=K(b)+n\mhf\xi_0+n\mo
\eta_0 ,
\end{equation}
where $\xi$ and $\eta$ are functions of two variables, $\xi_0$ and
$\eta_0$ are functions of a~single variable, each equals $O_P(1)$.
More specifically,
\begin{eqnarray*}
\xi(s,t)&=&{1\over n\half} \sumion(1-E) \{X_i(s)-\mu(s)\} \{
X_i(t)-\mu(t)\} ,\\
\xi_0(t)&=&{1\over n\half} \sumion(1-E) \{X_i(t)-\mu(t)\} \{
Y_i-E(Y_i)\} ,\\
\eta(s,t)&=&-n \{\bX(s)-\mu(s)\} \{\bX(t)-\mu(t)\} ,\\
\eta_0(t)&=&-n \{\bX(t)-\mu(t)\} (\bY-E\bY) .
\end{eqnarray*}
For any square-integrable function $L$ of two variables, define $\trip
L\trip^2=\intii L^2$ and put
$R_1=\trip K\trip+n\mhf\trip\xi\trip+n\mo\trip\eta\trip$,
$R_2=\trip\xi\trip+\trip\eta\trip$.
Define too
%
%e5.2 #&#
%
\begin{equation}\label{59}
\ze_j(t)=\inti K^j(b)(s) \xi(s,t) \,ds
\end{equation}
and
%
%e5.3 #&#
%e5.4 #&#
%
\begin{eqnarray}
\label{511}
\xi_j&=&K^{j-1}(\xi_0)+\sum_{k=0}^{j-2} K^k(\ze_{j-k-1}) ,
\\
\label{512}
\|\eta_j\|
&\leq& R_1^{j-1} \|\eta_0\|+R_2 \sum_{k=1}^{j-1} R_1^{j-k-1}
\bigl(\|K^k(b)\|+\bigtrip K^{k-1}\bigtrip\|\xi_0\|
\bigr)\nonumber\\[-8pt]\\[-8pt]
&&{}
+R_2 \trip\xi\trip\sum_{k=1}^{j-1} R_1^{j-k-1}
\sum_{\ell=0}^{k-2} \bigtrip K^\ell\bigtrip\|K^{k-\ell
-1}(b)\| .\nonumber
\end{eqnarray}

Theorem~\ref{Theorem3} below requires no assumptions beyond the model
at (\ref{21}), and the condition that
%
%e5.5 #&#
%
\begin{equation}\label{513}
\inti b^2<\infty,\qquad E\|X\|^4<\infty,\qquad E(\ep^2
)<\infty.
\end{equation}
[Recall that $\ep$, satisfying $E(\ep\mi X)=0$, is the error in the
model at (\ref{21}).] Note that, under (\ref{513}), it\vspace*{1pt} follows
from (\ref{511}) and (\ref{512}) that $\|\xi_j\|+\|\eta_j\|
=O_P(n\mhf)$. Theorem~\ref{Theorem3} shows that the empirical
approximations $\hkjb$ to the basis functions used by APLS, converge
in probability to their theoretical values $K^j(b)$ at a~rate $n\mhf$.
%
%th5.1 #&#
%
\begin{theo}\label{Theorem3}
If (\ref{513}) holds then, for each $j\geq1$,
%
%e5.6 #&#
%
\begin{equation}\label{510}
\hkjb=K^j(b)+n\mhf\xi_j+n\mo\eta_j ,
\end{equation}
where $\xi_j$ is defined at (\ref{511}) and $\eta_j$ satisfies
(\ref{512}).
\end{theo}

The next theorem shows that the matrix entries ${\hat h}_{jk}$, defined
at (\ref{52}), converge in probability to their theoretical
counterparts $h_{jk}$, at (\ref{54}), at a~rate $n\mhf$. This
theorem will be used to establish consistency of the empirical
coefficients $\hga_j$ used in the empirical APLS expansion at (\ref{35}).
Note that, since $\trip K\trip^2=\sumj\th_j^2$, the condition $0<\th
_1<\trip K\trip$ imposed in Theorem~\ref{Theorem4} is equivalent to
asserting that at least two values of $\th_j$ are nonzero. The
condition $\trip K\trip<1$ can be ensured by simply changing the scale
on which~$X$ is measured, and so is imposed without loss of generality.
%
%th5.2 #&#
%
\begin{theo}\label{Theorem4}
Assume (\ref{513}), that $\th_1,\th_2,\ldots$ is the eigenvalue
sequence in the representation (\ref{25}), ordered such that
(\ref{26}) holds, and that $0<\th_1<\trip K\trip<1$. Then $\|
\eta_j\|=O_p(\trip K\trip^j)$ uniformly in $1\leq j\leq C n\half$, and
%
%e5.7 #&#
%
\begin{eqnarray}\label{514}
{\hat h}_{jk}&=&h_{jk}+n\mhf\inti\{\xi_{j+1} K^k(b)+K^{j+1}(b)
\xi_k\}\nonumber\\[-8pt]\\[-8pt]
&&\hspace*{0pt}{}+O_p(n\mo\th_1^j \trip K\trip^k+n\mt\trip K\trip
^{j+k}),\nonumber
\end{eqnarray}
uniformly in $1\leq j\leq k\leq C n\half$ as $n\rai$, for each $C>0$.
\end{theo}

Our next result, Theorem~\ref{Theorem5}, applies Theorems~\ref{Theorem3} and~\ref{Theorem4} to derive a~stochastic expansion for the
difference between the theoretical approximant $g_p(x)$, at (\ref
{215}), and its estimator $\hg_p(x)$, at (\ref{35}).
Let $\De_{1jk}=\break\inti\{\xi_{j+1} K^k(b)+K^{j+1}(b) \xi_k\}$,
denoting the coefficient of $n\mhf$ in the expansion (\ref{514}),
and put $\De_1=(\De_{1jk})$, a~$p\times p$ matrix, and $\delta=(\De
_{101},\ldots,\De_{10p})\T$, a~$p$-vector. Also, let $\la=\la(p)$
be the smallest eigenvalue of the $p\times p$ matrix $H=(h_{jk})$,
introduced in Section~\ref{sec25}.
%
%th5.3 #&#
%
\begin{theo}\label{Theorem5}
Under the conditions of Theorem~\ref{Theorem4}, and if each $\th_j>0$,
%
%e5.8 #&#
%e5.9 #&#
%
\begin{eqnarray}
\label{515}
\|\hga-\{\ga+n\mhf H\mo(\delta-\De_1 \ga)\}
\|
&=&O_p(n\mo\la\mth) ,
\end{eqnarray}
\begin{eqnarray}\label{516}
&&\hg_p(x)-g_p(x)\nonumber\\
&&\qquad=\bY-EY
+n\mhf\sumjop\biggl[\{H\mo(\delta-\De_1 \ga)\}
_j\inti(x-EX)K^j(b)\nonumber\\[-8pt]\\[-8pt]
&&\qquad\quad\hspace*{95pt}{}+\ga_j\inti\{(x-EX) \xi_j-n\half(\bX-EX)
K^j(b)\}
\biggr]\nonumber\\
&&\qquad\quad{}+O_p(n\mo\la\mo\|\ga\|
+n^{-1} \la\mth) ,\nonumber
\end{eqnarray}
uniformly in functions $x$ and integers $p$ for which $\|x\|\leq C$,
$1\leq p\leq C n\half$ and $n\half\la\rai$, where $C>0$ is fixed
but arbitrary.
\end{theo}

Note that, by (\ref{54}), $|h_{jk}|\leq\th_1^{j+k+1} \|b\|^2$,
and therefore $\|Hv\|\leq C_1 \|v\|$ for all $p$-vectors $v$, where
the constant $C_1$ does not depend on $p$. (Here we have used the
condition $\th_1<1$, which we introduced in Theorem~\ref{Theorem4}
and also imposed in Theorem~\ref{Theorem5}.) Hence $\la\leq C_1$ for
all $p$. Note too that since, for finite $p$, $H$ is nonsingular (see
Section~\ref{sec25}), then its smallest eigenvalue $\lambda=\lambda
(p)$ is positive. On the other hand, when $p=\infty$ the sequence $\ga
_1,\ga_2,\ldots,$ that minimizes (\ref{213}) is not unique (see
Section~\ref{sec25}), and so we can have $\lambda\to0$ as $p\to
\infty$. The condition $n\half\la\rai$ imposed in Theorem~\ref{Theorem5} reflects this property, and essentially puts an upper bound
to the speed at which $p$ can tend to infinity as a~function of $n$.

%s5.3 #&#
\subsection{Implications of the main theorems and additional
results}\label{secfurther}

%s5.3.1 #&#
\subsubsection{Consistency and rates of convergence}
Let $X_0$ have the same distribution as $X_1,\ldots,X_n$ but be
independent of those random functions, and let \mbox{$\|\cdot\|\pred$}
denote the predictive $L_2$ norm, conditional on $X_1,\ldots,X_n$: if
$W$ is a~random variable, then
$
\|W\|\pred=\{E(W^2\bigmi X_1,\ldots,X_n)\}\half.
$
For example, taking $W=\hg_p(X_0)-g(X_0)$ we obtain a~measure of the
accuracy with which $\hg_p(X_0)$ predicts $g(X_0)$. We shall show in
Section~\ref{proof518} that if $p=p(n)$ is chosen to diverge no
faster than $n\half$, and sufficiently slowly to ensure that
%
%e5.10 #&#
%
\begin{equation}\label{517}
n\mhf\la\mo\|\ga\|+n\mo\la\mth\ra0
\end{equation}
as $n\rai$, then
%
%e5.11 #&#
%
\begin{eqnarray}\label{518}
&&\|\hg_p(X_0)-g(X_0)\|\pred\nonumber\\[-8pt]\\[-8pt]
&&\qquad=O_p\{n\mhf\la\mo(1+\|\ga\|)+n\mo\la\mth+t_p(\ga
_1,\ldots,\ga_p)\half\} ,\nonumber
\end{eqnarray}
where $t_p$ is as at (\ref{213}). It follows from Theorem~\ref{Theorem2} that if all of the eigenvalues $\th_j$ are nonzero, then
$t_p(\ga_1,\ldots,\ga_p)\ra0$ as $p\rai$. (As remarked in the
paragraph immediately below that theorem, the condition that each $\th
_j$ is nonzero can be dropped.) Therefore, (\ref{517}) implies that
$\hg_p(X_0)$ is consistent for $g(X_0)$.

Additionally,\vspace*{1pt} Theorems~\ref{Theorem3}--\ref{Theorem5} make it clear
that, provided $p$ does not diverge too quickly as a~function of $n$,
the quantities ${\sup_{j\leq p} }\|\hkjb-K^j(b)\|$, ${\sup_{1\leq
j,k\leq p}} |{\hat h}_{jk}-h_{jk}|$ and ${\sup_{j\leq p}} |\hga_j-\ga
_j|$ [see (\ref{521}) below] converge in probability to zero as $n$
diverges.

%s5.3.2 #&#
\subsubsection{Results in supremum metrics}
For our expansions of the function~$\hkjb$ at (\ref{510}), and of
the vector $\hga$ at (\ref{515}), our bounds on remainder terms
are given in $L_2$ metrics. In either case they can be extended to the
supremum metric. For example, (\ref{515}) itself implies that
%
%e5.12 #&#
%
\begin{equation}\label{521}
\sup_{1\leq j\leq p} |\hga_j-\ga_j
-n\mhf\{H\mo(\delta-\De_1 \ga)\}_j|
=O_p(n\mo\la\mth) .
\end{equation}

Theorem~\ref{Theorem6} below states a~version of (\ref{510}) in
the $L_\infty$ metric. It makes use of the following regularity conditions:
%
%e5.13 #&#
%
\begin{eqnarray}\label{522}
&&\mbox{for both $D_i\equiv1$ and $D_i\equiv Y_i$}\nonumber\\
&&\qquad\sup_{t\in\cI} \Biggl|{1\over n\half} \sumion\{X_i(t)
D_i-EX_i(t) D_i\}\Biggr|
=O_p(1),
\\
%
%e5.14 #&#
%
\label{523}
&&\sup_{t\in\cI} \inti\Biggl|{1\over n\half} \sumion(1-E) \{
X_i(s)-EX_i(s)\}
\{X_i(t)-EX_i(t)\}\Biggr|^2 \,ds\nonumber\\[-8pt]\\[-8pt]
&&\qquad=O_p(1) .\nonumber
\end{eqnarray}
Conditions (\ref{522}) and (\ref{523}) will be discussed in
Appendix~\ref{secA1}.
%
%th5.4 #&#
%
\begin{theo}\label{Theorem6}
If (\ref{513}), (\ref{522}) and (\ref{523}) hold, then
$\sup_{t\in\cI} |\xi_j(t)|=O_p(1)$ for each $j$, and
\[
\sup_{t\in\cI} |\hkjb(t)-\{K^j(b)(t)+n\mhf\xi
_j(t)\}|
=O_p(n\mo) .
\]
\end{theo}

%s5.3.3 #&#
\subsubsection{Interpreting stochastic expansions}
The coefficients\vspace*{1pt} of $n\mhf$ in the expansions of $\hkjb
(t)-K^j(b)(t)$, ${\hat h}_{jk}-h_{jk}$, $\hga_j-\ga_j$ and $\hga
_p(x)-\ga_p(x)$ in~(\ref{510}), (\ref{514}), (\ref{515})
[see also (\ref{521})] and (\ref{516}), respectively, are each
equal to~$n\mo$ multiplied by a~sum of $n$ independent and identically
distributed random variables with zero mean, plus a~term that equals
$O_p(n\mo)$. In these cases, for fixed $(j,t)$, $(j,k)$, $j$ and
$(p,x)$, respectively, the independent random variables do not depend
on $n$. Therefore, their variances can be computed easily.

For example, in the case of ${\hat h}_{jk}-h_{jk}$, using (\ref{514})
and the definitions of $\xi$ and~$\xi_0$, we have, under the
conditions of Theorem~\ref{Theorem4} and for each fixed $j$ and $k$,
%
%e5.15 #&#
%
\begin{equation}\label{524}
{\hat h}_{jk}=h_{jk}+n\mo\sumion Z_{ijk}
+O_p(n\mo) ,
\end{equation}
where the independent and identically distributed random variables
$Z_{1jk},\ldots,\break Z_{njk}$ are given by
\begin{eqnarray*}
Z_{ijk}
&=&(1-E)\\
&&{}\times\inti\Biggl(K^k(b)(u)
\Biggl[\{Y_i-E(Y_i)\}\\
&&\hspace*{-24.6pt}\qquad\hspace*{79.1pt}{}\times \inti K^j(u,t) \{X_i(t)-\mu(t)\} \,dt\\
&&\hspace*{-24.6pt}\qquad\hspace*{79.1pt}{}
+\sum_{\ell=0}^{j-1} \inti\{X_i(t)-\mu(t)\} K^\ell(t,u) \,dt\\
&&\hspace*{-24.6pt}\qquad\hspace*{106.2pt}{}\times
\inti K^{j-\ell}(b)(s) \{X_i(s)-\mu(s)\} \,ds\Biggr]\\
&&\hspace*{-24.6pt}\hspace*{53.3pt}{}
+K^{j+1}(b)(u) \Biggl[\{Y_i-E(Y_i)\}\\
&&\hspace*{-24.6pt}\hspace*{124pt}{}\times
\inti K^{k-1}(u,t) \{X_i(t)-\mu(t)\} \,dt\\
&&\hspace*{-24.6pt}\hspace*{124pt}{}
+\sum_{\ell=0}^{k-2} \inti\{X_i(t)-\mu(t)\} K^\ell(t,u) \,dt\\
&&\hspace*{-24.6pt}\qquad\qquad\hspace*{106pt}{}\times
\inti K^{k-\ell-1}(b)(s) \{X_i(s)-\mu(s)\} \,ds\Biggr]\Biggr) \,du .
\end{eqnarray*}
The distribution of $Z_{ijk}$ does not depend on $n$, and, under the
assumption of finite fourth moment of $X$ and finite second moment of
$\ep$ [see (\ref{513})], $Z_{ijk}$ has finite variance $\si_{jk}
^2$, say. Hence, for each fixed $j$ and $k$ it follows from~(\ref
{524}) that $n\half({\hat h}_{jk}-h_{jk})$ is asymptotically normal
N$(0,\si_{jk}^2)$.

%s5.3.4 #&#
\subsubsection{Hankel matrix properties} In Section~\ref{sec25} we
demonstrated that $\a_j=\int x^j m(dx)$, where $m$ is the measure
that places mass $(\be_r\PC)^2 \th_r$ at the point $\th_r$ for
$r\geq1$; $m$ has no mass anywhere else. Therefore the $p\times p$
matrix $H=(\a_{j+k})$ is a~Hankel matrix for which the associated
nonnegative measure, $m$, is discrete and compactly supported. The
latter property implies that $m$ is completely determined by its
moments $\a_j$, and hence that the Hankel matrix is ``determinate;''
see, for example, \citet{BerSzw11}. In such cases the smallest
eigenvalue of $H$ can converge to zero arbitrarily fast as $p$ diverges
[\citet{BerSzw11}, Theorem~2.5], although more is known about the
case where $m$ is a~continuous than that of a~discrete measure, and it
is particularly challenging to develop general theory describing
properties of $H\mo$ in the context of our measures $m$. [See
\citet{Las90} and \citet{HouLasMu05} for access to the
literature on inverses of
Hankel matrices and their determinants.] Nevertheless, as we noted in
Section~\ref{sec25}, $H$ is generally nonsingular for all $p$.

%s6 #&#
\section{Numerical illustrations}\label{sec4}
In this section we illustrate, numerically, in a~few examples, the fact
that the algorithms in Section~\ref{secaltAPLS} and Appendix~\ref{sec31} do indeed solve the same problem. We also illustrate the main
difference between the PLS basis and the PCA basis, namely that PLS can
capture the interaction between $X$ and $Y$ using a~smaller number of
terms than PCA.

In our first illustration, we take the $X_i$'s from a~real data study,
and generate the $Y_i$'s according to the linear model at (\ref
{21}). By choosing the population in this way, we can represent, in
simulations, the vagaries of real data, but we can still compare the
performance of our methodology with the ``truth.'' We take the $X_i$
curves from a~benchmark Phoneme dataset, which can be downloaded from
\href{http://www-stat.stanford.edu/\textasciitilde tibs/ElemStatLearn/}{www-stat.stanford.edu/ElemStatLearn}. In these
data, $X_i(t)$
represents log-periodograms constructed from recordings of different
phonemes. The periodograms are available at 256 equispaced frequencies
$t$, which for simplicity we denote by $t=1,2,\ldots,256$. Hence, in
this example, $\cI=[1,256]$. See \citet{HasTibFri09} for more
information about this dataset. We used the $N=1717$ data
curves~$X_i(t)$ that correspond to the phonemes ``aa'' as in ``dark'' and
``ao'' as in ``water.''\looseness=1

We computed the first $J=20$ empirical PCA basis functions $\hat\phi
_1(t),\ldots,\break\hat\phi_{20}(t)$, and considered four different
curves $b$, which we constructed by taking $b(t)=\sum_{j=1}^{J}a_j\hat
\phi_j(t)$ for four different sequences of $a_j$'s: (i)
$a_j=(-1)^j\cdot\break1\{j\leq5\}$; (ii)~$a_j=(-1)^j\cdot1\{6\leq j\leq
10\}$; (iii) $a_j=(-1)^j\cdot1\{11\leq j\leq15\}$;
(iv)~$a_j=(-1)^j\cdot1\{16\leq j\leq20\}$.
These four models were chosen to illustrate clearly the advantages of
the PLS basis over the PCA basis. Example (i) illustrates a~situation
particularly favourable to PCA, where the interaction between $X$ and
$Y$ can be represented by the first few PCA basis functions. There we
do not expect that PCA will need many more terms than PLS to achieve a~small prediction error. On going from example (i) to example (iv), the
function $b$ is represented by five consecutively indexed PCA basis
functions in each case, but with their indices successively larger.
However, as we shall see below, in those cases too, PLS manages to
construct a~basis that captures the interaction between $X$ and $Y$
using only the first few terms.

In the four cases, for $i=1,\ldots,N$ we generated the $Y_i$'s by
taking $Y_i=\inti X_i b+\varepsilon_i$, where $\varepsilon_i\sim
\mathrm{N}(0,\sigma^2)$, and where $5\sigma^2$ was equal to the empirical
variance of the $\inti bX_i$'s calculated from the $N$ observations.
Then, in each case, we randomly split these $N$ observations in two
parts: a~training sample of size $n$, and a~test sample of size $N-n$.
We did this 200 times for each of $n=30$, $n=50$ and $n=100$, so for
each setting we generated $200$ test and training samples.

For each set of test and training samples generated in this way, we
constructed our predictor using only the test sample, and
then\vspace*{1pt} we applied it to predict $\inti b X_i$ for each $X_i$
in the associated training sample. In other words, we constructed $\bar
X$, $\bar Y$ and $\hat b$ from the training sample only, where $\hat b$
was the empirical version of $b_p$ calculated either via the first $p$
terms of the PLS basis (calculated from the algorithm in Appendix
\ref{sec31} or the second algorithm of Section~\ref{secaltAPLS}), or
via the first $p$ terms of the PCA basis, for each of $p=1,\ldots,10$.
Then, for each observation $X_i$ in the test sample, we calculated the
predictor $\hat Y_i=\bar Y+\inti\hat b(X_i-\bar X)$ of $\inti b X_i$.
Note that this predictor includes the estimator $\bar Y-\inti b\bar
X_i$ of the intercept because, although our data were generated from a~model with no intercept, in practice we are not supposed to know this.

%f1 #&#
%
\begin{figure}

\includegraphics{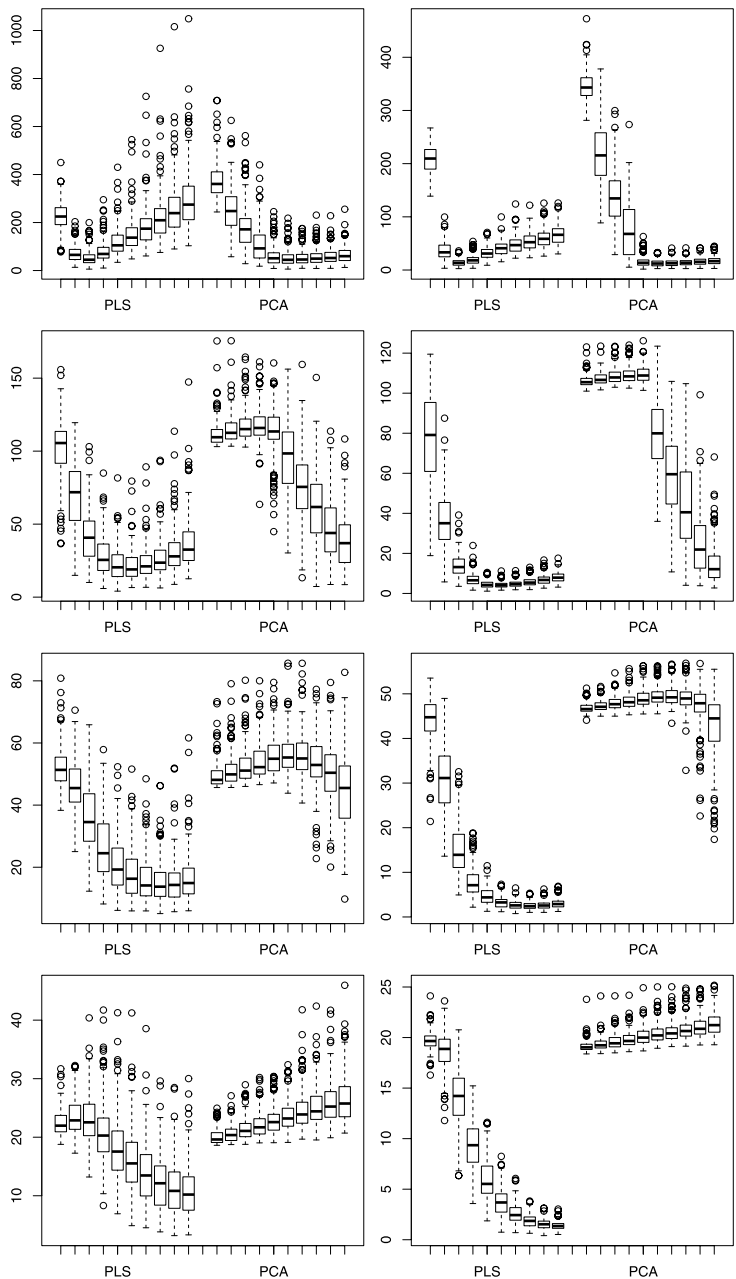}

\caption{Boxplots of the prediction error using the
first $p$ PLS components (first group of 10 boxes) or the first $p$ PCA
components (last group of 10 boxes), calculated from 200 samples of
sizes $n=30$ (first column) or $n=100$ (second column) generated from
the phoneme data. The curve $b$ is that in cases \textup{(i), (ii),
(iii)} and \textup{(iv)}, in, respectively, rows 1, 2, 3 and 4. From
left to right, each group of 10 boxplots addresses the settings indexed
by $p=1$ to $p=10$.}\label{figPhoneme}
\end{figure}

\begin{figure}

\includegraphics{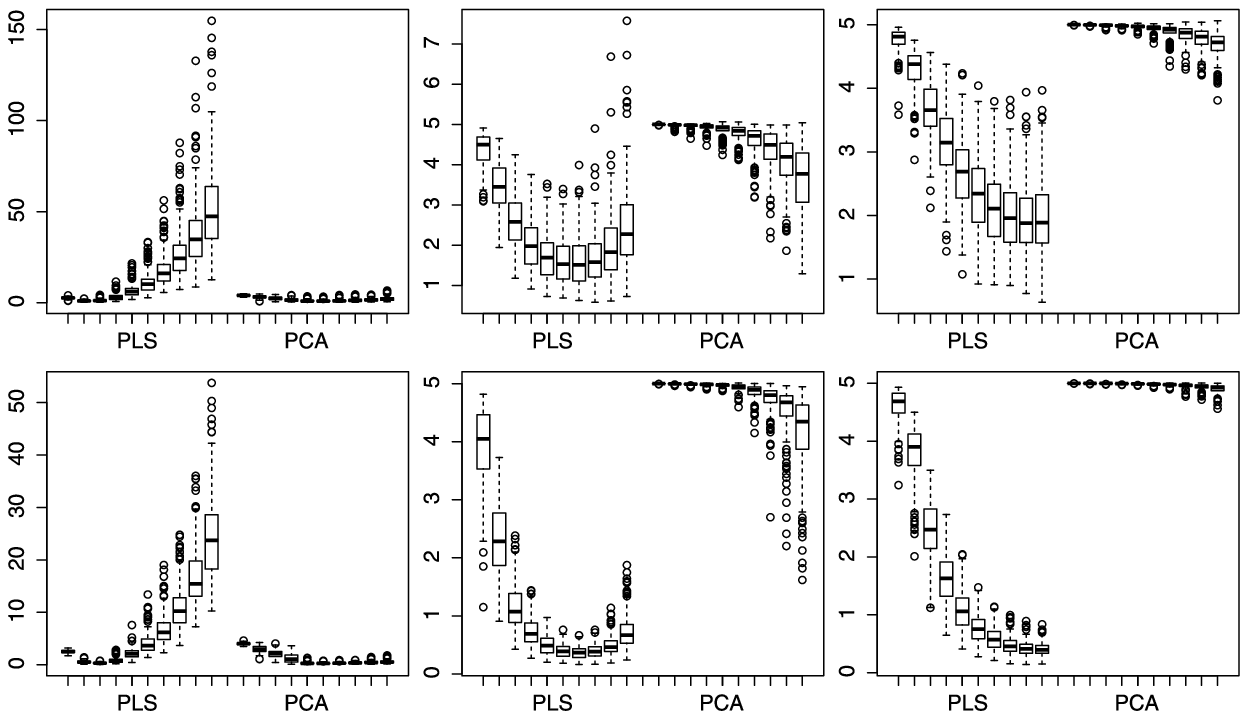}

\caption{Boxplots of the ISE of $\hat b$ using the first $p$ PLS
components (first group of 10 boxes) or the first $p$ PCA components
(last group of 10 boxes), calculated from 200 samples of sizes $n=30$
(first row) or $n=100$ (second row) generated from the phoneme data.
The curve $b$ is that in cases \textup{(i), (iii)} and \textup{(iv)},
in, respectively,
columns 1, 2 and 3. From left to right, each group of 10 boxplots
addresses the settings indexed by $p=1$ to
$p=10$.}\label{figPhonemeb}
\end{figure}

To quantify\vspace*{1pt} the quality of prediction, we calculated the prediction
error $\mathrm{PE}=(N-n)^{-1}\sum_{i=1}^{N-n}(\hat Y_i-\inti b
X_i)^2$ in each case, for each method, and for each test sample. In
Figure~\ref{figPhoneme} we show boxplots of these prediction errors
calculated in each case from the $200$ test samples. Note that here the
two PLS algorithms gave exactly the same estimators, and so the
boxplots only show the results for the standard PLS algorithm and for
the PCA method.
These boxplots show that as the information about the interaction
between~$X$ and~$Y$ moves further away in the sequence of $\hat\phi
_j$'s [i.e., going from case (i) to case (iv)], PLS can capture the
interaction using fewer terms than PCA. For example, in case (i), PLS
took $p=3$ components to reach the prediction error that PCA reached
with $p=5$, but in case (iv), the prediction error was already very
small for PLS with $p=10$, and was still very large\vspace*{1pt} for PCA with $p=10$.
We also calculated the integrated squared error $\mathrm{ISE}=\inti
(\hat b -b )^2$ for each method and test sample. In Figure~\ref{figPhonemeb} we show boxplots of these ISEs calculated from the $200$
test samples, for models (i), (iii) and (iv). We can see that the PLS
estimator of $b$ needs fewer components than the PCA estimator to reach
small ISE values.

%f2 #&#
%

%f3 #&#
%
\begin{figure}[b]

\includegraphics{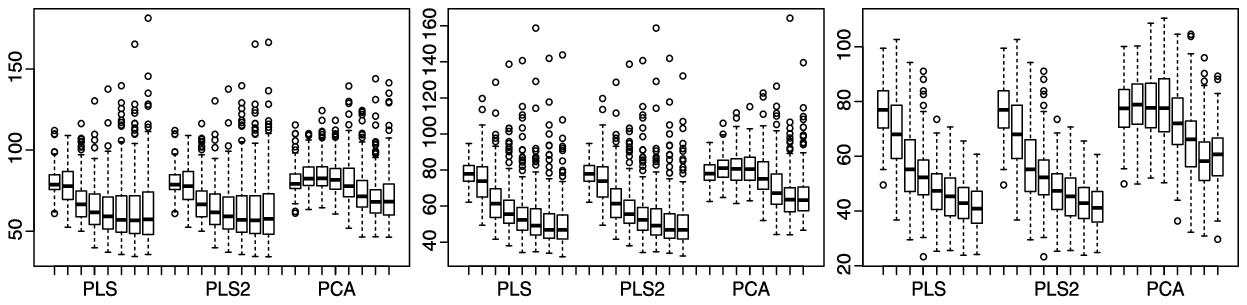}

\caption{Boxplots of the estimated prediction error using the first $p$
PLS components calculated by the algorithm of Appendix \protect\ref{sec31}
(first group of 8 boxes) or the second algorithm of Section
\protect\ref{secaltAPLS} (second group of 8 boxes, denoted by PLS2),
or the
first $p$ PCA components (last group of 8 boxplots). Each box was
calculated from 200 samples of sizes $n=30$ (first column), $n=50$
(second column) or $n=100$ (third column) drawn randomly from the
orange data. From left to right, each group of 8~boxplots is for $p=1$
to $p=8$.}\label{figOrange}
\end{figure}

In our second example we took the orange juice data which comprise
$N=216$ observations $(X_i(t),Y_i)$, $i=1,\ldots,N$, where
each $Y_i$ is the saccharose content of a~sample of orange juice, and
$X_i$ is a~curve representing the first derivative of near-infrared
spectra of the juice at 700 equispaced points $t$. We take $t=1,\ldots
,700$ (hence $\cI=[1,700]$). The data can be found at
\href
{http://www.ucl.ac.be/mlg/index.php?page=DataBases}{www.ucl.ac.be/mlg/index.php?\allowbreak page=DataBases}.
As with our simulated
data above, we split the observations randomly into a~training sample
of size~$n$ and a~test sample of size $N-n$, for each of $n=30$, $50$
and $100$. We did this 200 times for each $n$. Then in each case we
calculated our predictor, as above, from the training sample, and
applied it for predicting $\inti X_i b$ for the corresponding test
sample. Here we did not know the true model, so we calculated an
estimator of the prediction error as $\hat{\mathrm{PE}}
=(N-n)^{-1}\sum_{i=1}^{N-n}(\hat Y_i-Y_i)^2$, for each
$(X_i,Y_i)$ in the test sample. In this way we obtained 200 values of
$\hat{\mathrm{PE}}$ for each $n$. Figure~\ref{figOrange} shows, for
each $n$, boxplots of these 200 $\hat{\mathrm{PE}}$'s, for $p=1$ to
$8$. As above, the two PLS algorithms (the algorithm in Appendix~\ref{sec31} and the second algorithm of Section~\ref{secaltAPLS}) gave
exactly the same results, except for $p=8$ where the numerical
roundings of both methods differed somewhat. Therefore we show the
boxplots for both algorithms. In this example too we can see that the
two PLS algorithms clearly solve the same problem, and that PLS needs
fewer terms (i.e., $p$ is smaller) to capture the same interactions as
PCA. This can be advantageous when computing time is an issue, for
example when a~linear prediction is associated with a~complex
nonparametric procedure. For example, in
\citet{FerVie06}, the
linear fit is used in combination with nonparametric estimators of $E(Y|X)$.

%s7 #&#
\section{Technical arguments}\label{sec6}

%s7.1 #&#
\subsection{\texorpdfstring{Proof of Theorem \protect\ref{Theorem1}}{Proof of Theorem 3.1}}\label{sec61}
Defining $\si^2=\var(\ep)$ we see that the right-hand side of (\ref
{27}) can be expressed as
\begin{eqnarray*}
&&\cov\Biggl\{\biggl(\inti b X\biggr)
-\sum_{j=1}^{p-1} \biggl(\inti b \psi_j\biggr) \biggl(\inti X
\psi_j\biggr),
\inti X \psi_p\Biggr\}\\
&&\qquad=\intii b \psi_p K
-\sum_{j=1}^{p-1} \biggl(\inti b \psi_j\biggr)
\biggl(\intii\psi_j \psi_p K\biggr) .
\end{eqnarray*}
The partial derivative of the right-hand side here, with respect to
$\psi_p$, equals
%
%e7.1 #&#
%
\begin{equation}\label{61}
K\Biggl\{b-\sum_{j=1}^{p-1}
\biggl(\inti b \psi_j\biggr) \psi_j\Biggr\}.
\end{equation}
The equation in $c_k$ at (\ref{210b}) is the result of adjoining
Lagrange multipliers on the right-hand side so as to accommodate the
first $p-1$ constraints in~(\ref{28b}). The factor $c_0$ on the
right-hand side of (\ref{29}) accommodates the last constraint
in (\ref{28b}).

%s7.2 #&#
\subsection{\texorpdfstring{Proof of Theorem \protect\ref{Theorem2}}{Proof of Theorem 3.2}}\label{sec62}
Recall that $\cC(\cI)$ is the space of all square-integrable
functions on $\cI$, and suppose $b=\sumj\be\PC_j \phi_j\in\cC
(\cI)$. Write $\cC_p(\cI)$ for the $p$-dimen\-sional subspace of
$\cC(\cI)$ generated by the PCA basis func\-tions~$\phi_1,\ldots,\phi
_p$, and let $K_p$ denote the transformation that takes $b_p\equiv\break\sum
_{1\leq j\leq p} \be\PC_j \phi_j\in\cC_p(\cI)$ to $\sum_{1\leq
j\leq p} \th_j\be\PC_j \phi_j$.
Now,
\[
(\th_1 I-K_p)\cdots(\th_p I-K_p) b_p=0
\]
for all $b_p\in\cC_p(\cI)$. Therefore,
%
%e7.2 #&#
%
\begin{equation}\label{62}
a_0 b_p+a_1 K_p(b_p)+\cdots+a_p K_p^p(b_p)=0
\end{equation}
for all $b_p\in\cC_p(\cI)$, where $a_0,\ldots,a_p$ are constants
and $a_0=\th_1\cdots\th_p$. In particular, $a_0$ is nonzero, and so
(\ref{62}) implies that, for constants $c_1,\ldots,c_p$,
%
%e7.3 #&#
%
\begin{equation}\label{63}
b_p=c_1 K(b_p)+\cdots+c_p K^p(b_p) .
\end{equation}

Let $P_p$ denote the projection operator that takes $b=\sumj\be\PC
_j \phi_j\in\cC(\cI)$ to $P_p(b)=b_p\in\cC_p(\cI)$. Since $P_p$
and $K$ commute, then $K^j(b_p)=K^j P_p(b)=P_p K^j(b)$.
Therefore (\ref{63}) implies that $b_p=P_p\{c_1 K(b)+\cdots+c_p
K^p(b)\}$, or equivalently,
%
%e7.4 #&#
%
\begin{equation}\label{64}
P_p[b-\{c_1 K(b)+\cdots+c_p K^p(b)\}]=0 .
\end{equation}
In view of (\ref{64}), if we let $\cD(\cI)$ denote the vector
space generated by $K(b)$, $K^2(b),\ldots,$ and if we define $P_p\{
\cD(\cI)\}=\{P_p(z)\dvtx z\in\cD(\cI)\}$, then $P_p(b)\in P_p\{\cD(\cI
)\}$ for all\vadjust{\goodbreak} $p$. Now, $P_p\{\cD(\cI)\}\subseteq\cD(\cI)$, which
is closed under limit operations in $L_2$. Therefore, the limit as
$p\to\infty$ of $P_p(b)$, that is $b$, must be in $\cD(\cI)$.

%s7.3 #&#
\subsection{\texorpdfstring{Proof of Theorem \protect\ref{Theorem3}}{Proof of Theorem 5.1}}\label{sec63}
Assume it can be proved that (\ref{510}) holds, with $\xi_j$ and
$\eta_j$ satisfying (\ref{511}) and (\ref{512}), for a~particular $j\geq1$; in view of~(\ref{58}), (\ref{510}) is
valid for $j=1$. Then,
%
%e7.5 #&#
%
\begin{eqnarray}\label{65}
\hkjob(t)&=&\inti\hkjb(s) \hK(s,t) \,ds\nonumber\\[-2pt]
&=&\inti\{K^j(b)+n\mhf\xi_j+n\mo\eta_j\}(s)(K+n\mhf
\xi+n\mo\eta)(s,t) \,ds\nonumber\\[-2pt]
&=&K^{j+1}(b)(t)+n\mhf\inti\{K^j(b)(s) \xi(s,t)+\xi_j(s)
K(s,t)\} \,ds\nonumber\\[-2pt]
&&{}
+n\mo\inti\{K^j(b)(s) \eta(s,t)+\eta_j(s) K(s,t)+\xi_j(s)
\xi(s,t)\} \,ds\nonumber\\[-2pt]
&&{}
+n^{-3/2}\inti\{\xi_j(s) \eta(s,t)+\eta_j(s) \xi(s,t)\}
\,ds\nonumber\\[-2pt]
&&{}
+n\mt\inti\eta_j(s) \eta(s,t) \,ds .
\end{eqnarray}
Therefore, taking $\xi_{j+1}$ to be given by the coefficient of $n\mhf
$ in (\ref{65}), and recalling the definition of $\ze_j$ at (\ref
{59}), we have
\begin{eqnarray*}
\xi_{j+1}(t)&=&\inti\{K^j(b)(s) \xi(s,t)+\xi_j(s) K(s,t)
\} \,ds\\[-2pt]
&=&K(\xi_j)(t)+\ze_j(t)
=K^2(\xi_{j-1})(t)+K(\ze_{j-1})(t)+\ze_j(t) ,
\end{eqnarray*}
which, on iteration, gives (\ref{511}).

Finally we derive the bound at (\ref{512}) on the remainder, again
arguing by induction; assuming that (\ref{512}) holds for $j$ we
establish it for $j+1$. Taking~$\eta_{j+1}$ to equal $n$ times the sum
of the terms in $n\mo$, $n^{-3/2}$ and $n\mt$ in~(\ref{65}), we
deduce that
\begin{eqnarray*}
\eta_{j+1}(t)&=&\inti\{K^j(b)(s) \eta(s,t)+\eta_j(s)
K(s,t)+\xi_j(s) \xi(s,t)\} \,ds\\[-2pt]
&&{}
+n\mhf\inti\{\xi_j(s) \eta(s,t)+\eta_j(s) \xi(s,t)\} \,ds
+n\mo\inti\eta_j(s) \eta(s,t) \,ds .
\end{eqnarray*}
Therefore,
\begin{eqnarray*}
\|\eta_{j+1}\|
&\leq&\|K^j(b)\| \trip\eta\trip+\|\eta_j\| \trip K\trip
+\|\xi_j\| \trip\xi\trip%\\
%&
+n\mhf(\|\xi_j\| \trip\eta\trip+\|\eta_j\| \trip\xi\trip)\\
&&{} +n\mo\|\eta_j\| \trip\eta\trip\\
&\leq&\bigl(\|K^j(b)\|+\|\xi_j\|\bigr)R_2
+\|\eta_j\| R_1\\
&\leq&\bigl(\|K^j(b)\|+\|\xi_j\|\bigr)R_2
+\bigl\{\bigl(\|K^{j-1}(b)\|+\|\xi_{j-1}\|\bigr)R_2
+\|\eta_{j-1}\|R_1\bigr\}R_1\\
&=&\bigl\{\bigl(\|K^j(b)\|+\|\xi_j\|\bigr)
+\bigl(\|K^{j-1}(b)\|+\|\xi_{j-1}\|\bigr) R_1\bigr\} R_2
+\|\eta_{j-1}\| R_1^2\\
&\leq&\Biggl\{\sum_{k=1}^j \bigl(\|K^k(b)\|+\|\xi_k\|
\bigr) R_1^{j-k}\Biggr\} R_2
+\|\eta_1\| R_1^j ,
\end{eqnarray*}
where the last identity follows on iteration. Observe too that, by
(\ref{58}), $\eta_1=\eta_0$. Therefore,
%
%e7.6 #&#
%
\begin{equation}\label{66}
\|\eta_{j+1}\|
\leq R_1^j \|\eta_0\|
+R_2 \sum_{k=1}^j R_1^{j-k} \bigl(\|K^k(b)\|+\|\xi_k\|\bigr) .
\end{equation}

Note too that, by (\ref{59}), $\|\ze_j\|\leq\|K^j(b)\| \trip\xi
\trip$, and so, by (\ref{511}),
\[
\|\xi_k\|\leq\bigtrip K^{k-1}\bigtrip\|\xi_0\|
+\trip\xi\trip\sum_{\ell=0}^{k-2} \bigtrip K^\ell\bigtrip
\|K^{k-\ell-1}(b)\| .
\]
Hence, by (\ref{66}),
%
%e7.7 #&#
%
\begin{eqnarray}\label{67}
\|\eta_{j+1}\|-R_1^j \|\eta_0\|
&\leq& R_2 \sum_{k=1}^j R_1^{j-k}
\Biggl\{\|K^k(b)\|
+\bigtrip K^{k-1}\bigtrip\|\xi_0\|\nonumber\\
&&\hspace*{61pt}{}
+\trip\xi\trip\sum_{\ell=0}^{k-2} \bigtrip K^\ell\bigtrip
\|K^{k-\ell-1}(b)\|\Biggr\}\nonumber\\
&=&R_2 \sum_{k=1}^j R_1^{j-k}
\bigl(\|K^k(b)\|+\bigtrip K^{k-1}\bigtrip\|\xi_0\|
\bigr)\nonumber\\
&&{}
+R_2 \trip\xi\trip\sum_{k=1}^j R_1^{j-k}
\sum_{\ell=0}^{k-2} \bigtrip K^\ell\bigtrip\|K^{k-\ell
-1}(b)\| .
\end{eqnarray}
Result (\ref{512}) for $j+1$ follows from (\ref{67}).

%s7.4 #&#
\subsection{\texorpdfstring{Proof of Theorem \protect\ref{Theorem4}}{Proof of Theorem 5.2}}\label{sec64}
Representation (\ref{510}) implies that
%
%e7.8 #&#
%
\begin{equation}\label{68}
{\hat h}_{jk}=h_{jk}
+n\mhf\inti\{\xi_{j+1} K^k(b)+K^{j+1}(b) \xi_k\}
+n\mo R_{jk} ,
\end{equation}
where
%
%e7.9 #&#
%
\begin{eqnarray}\label{69}
{\hat h}_{jk}&=&\inti\hkjob\hkkb,\qquad h_{jk}=\inti K^{j+1}(b)
K^k(b) ,\nonumber\\
|R_{jk}|
&\leq& \biggl|\inti\{\xi_{j+1} \xi_k+K^{j+1}(b) \eta
_k+K^k(b) \eta_{j+1}\}\nonumber\\
&&{}
+n\mhf\inti(\eta_{j+1} \xi_k+\xi_{j+1} \eta_k)
+n\mo\inti\eta_{j+1} \eta_k\biggr|\nonumber\\
&\leq&\|\xi_{j+1}\| \|\xi_k\|+\|K^{j+1}(b)\| \|\eta_k\|
+\|K^k(b)\| \|\eta_{j+1}\|\nonumber\\
&&{}
+n\mhf(\|\eta_{j+1}\| \|\xi_k\|+\|\xi_{j+1}\| \|\eta_k\|)
+n\mo\|\eta_{j+1}\| \|\eta_k\| .
\end{eqnarray}

Next we bound $|R_{jk}|$. Observe that $\trip K^k\trip^2=\sumj\th
_j^{2k}=O(\th_1^{2k})$, $\|K^k(b)\|^2=\sumj\th_j^{2k} (\inti
b\phi_j)^2=O(\th_1^{2k})$ and $\trip\eta\trip+\trip\xi\trip
+\|\eta_0\|+\|\xi_0\|=O_p(1)$ as $n\rai$. Hence, by (\ref{59}),
$\|\ze_j\|\leq\|K^j(b)\| \trip\xi\trip=O_p(\th_1^j)$, uniformly
in $j\geq1$, and therefore by (\ref{511}),
%
%e7.10 #&#
%
\begin{equation}\label{610}
\|\xi_j\|=O_p\Biggl(\th_1^j+\sum_{k=0}^{j-2} \th_1^k \th
_1^{j-k-1}\Biggr)
=O_p(j \th_1^j) ,
\end{equation}
uniformly in $j\geq1$. Note too that
\[
R_1^j=\{(\trip K\trip+n\mhf\trip\xi\trip+n\mo\trip
\eta\trip)^j\}
=O_p(\trip K\trip^j) ,
\]
uniformly in $1\leq j\leq C n\half$, for any $C>0$. More simply,
$R_2=O_p(1)$. Hence, by (\ref{512}),
%
%e7.11 #&#
%
\begin{eqnarray}\label{611}
\|\eta_j\|
&=&O_p\Biggl(\trip K\trip^j+\sum_{k=1}^{j-1} \trip K\trip^{j-k-1}
\th_1^k
+\sum_{k=1}^{j-1} \trip K\trip^{j-k-1} k
\th_1^k\Biggr)\nonumber\\
&=&O_p(\trip K\trip^j) ,
\end{eqnarray}
uniformly in $1\leq j\leq C n\half$. (Here we have used the property
$0<\th_1<\trip K\trip<1$.) Combining (\ref{69})--(\ref{611})
we find that
%
%e7.12 #&#
%e7.13 #&#
%
\begin{eqnarray}
R_{jk}&=&O_p\{jk \th_1^{j+k}+\th_1^j \trip K\trip^k+\th_1^k
\trip K\trip^j\\
&&\hspace*{16.5pt}{}
+n\mhf(\trip K\trip^j \th_1^k+\trip K\trip^k \th_1^j)
+n\mo\trip K\trip^{j+k}\}\nonumber\\
\label{612}
&=&O_p(\th_1^j \trip K\trip^k+n\mo\trip K\trip^{j+k}),
\end{eqnarray}
uniformly in $1\leq j\leq k\leq C n\half$. Theorem~\ref{Theorem4}
follows from (\ref{68}) and (\ref{612}).

%s7.5 #&#
\subsection{\texorpdfstring{Proof of Theorem \protect\ref{Theorem5}}{Proof of Theorem 5.3}}\label{sec65}
From (\ref{215}), (\ref{216}) and (\ref{35}) we deduce that
%
%e7.14 #&#
%
\begin{eqnarray}\label{613}
&&\hg_p(x)-g_p(x)-(\bY-EY)\nonumber\\
&&\qquad=\sumjop\biggl\{\hga_j \inti(x-\bX) \hkjb
-\ga_j \inti(x-EX) K^j(b)\biggr\}\nonumber\\
&&\qquad=\sumjop\biggl[(\hga_j-\ga_j) \inti(x-EX) K^j(b)\nonumber\\
&&\qquad\quad\hspace*{19.5pt}{}+\ga_j\inti(x-EX)
\{\hkjb-K^j(b)\}\nonumber\\
&&\qquad\quad\hspace*{19.5pt}{}-\ga_j\inti(\bX-EX) K^j(b)\nonumber\\
&&\qquad\quad\hspace*{19.5pt}{}+(\hga_j-\ga_j)\inti(x-EX) \{\hkjb-K^j(b)\}
\nonumber\\
&&\qquad\quad\hspace*{19.5pt}{}
-(\hga_j-\ga_j)\inti(\bX-EX) K^j(b)\nonumber\\
&&\qquad\quad\hspace*{49pt}{}-\hga_j\inti(\bX-EX) \{\hkjb-K^j(b)\}\biggr].
\end{eqnarray}

Combining\vspace*{1pt} (\ref{510}), (\ref{610}) and the bound $\|\eta_j\|
=O_p(\trip K\trip^j)$, valid uniformly in $1\leq j\leq C n\half$ for
each $C>0$ and given in Theorem~\ref{Theorem4}, we deduce that
%
%e7.15 #&#
%
\begin{eqnarray}\label{614}
\|\hkjb-K^j(b)\|
&\leq& n\mhf\|\xi_j\|+n\mo\|\eta_j\|\nonumber\\[-8pt]\\[-8pt]
&=&O_p(n\mhf j \th_1^j
+n\mo\trip K\trip^j) ,\nonumber
\end{eqnarray}
uniformly in $1\leq j\leq C n\half$ for each $C>0$.

More simply, $\|\bX-EX\|=O_p(n\mhf)$. Combining this bound with
(\ref{610}), (\ref{613}) and the properties $\|x\|\leq C$ and
$\|K^j(b)\|=O(\th_1^j)$, we deduce that
%
%e7.16 #&#
%
\begin{eqnarray}\label{615}
&&\hg_p(x)-g_p(x)-(\bY-EY)\nonumber\\
&&\qquad=\sumjop\biggl[(\hga_j-\ga_j) \inti(x-EX)
K^j(b)\nonumber\\
&&\qquad\quad\hspace*{19pt}{}
+\ga_j\inti(x-EX) \{\hkjb-K^j(b)\}
-\ga_j\inti(\bX-EX) K^j(b)\biggr]\nonumber\\
&&\qquad\quad{}
+O_p\Biggl\{n\mhf\sumjop
(|\hga_j-\ga_j|+n\mhf|\hga_j|)
(j \th_1^j+n\mhf\trip K\trip^j)\Biggr\} ,
\end{eqnarray}
uniformly in $1\leq p\leq C n\half$ and $\|x\|\leq C$, for each
$C>0$. Using (\ref{510}) and the bound $\|\eta_j\|=O_p(\trip
K\trip^j)$, we deduce from (\ref{615}) that
%
%e7.17 #&#
%
\begin{eqnarray}\label{616}
&&\hg_p(x)-g_p(x)-(\bY-EY)\nonumber\\
&&\qquad=\sumjop\biggl[(\hga_j-\ga_j) \inti(x-EX)
K^j(b)\nonumber\\
&&\qquad\quad\hspace*{19pt}{} +\ga_j\inti\{n\mhf(x-EX) \xi_j-(\bX-EX)
K^j(b)\}\biggr]\nonumber\\
&&\qquad\quad{} +O_p\Biggl[n\mhf\sumjop
\{|\hga_j-\ga_j| (j \th_1^j+n\mhf\trip K\trip^j
)\nonumber\\
&&\qquad\quad\hspace*{90pt}{}
+n\mhf(|\hga_j|+|\ga_j|) \trip K\trip^j\}\Biggr] ,
\end{eqnarray}
uniformly in $1\leq p\leq C n\half$ and $\|x\|\leq C$, for each $C>0$.

Given any $p\times p$ matrix $M$, define its norm by $\|M\|=\sup_{v
\dvtx\|v\|=1} \|Mv\|$. Writing $\De$ for a~particular $p\times p$
matrix, and recalling that $\la=\la(p)$ denotes the smallest eigenvalue
of $H$, we have $\|\De H\mo\|\leq\|\De\|/\la$. Therefore, if
$\hH=({\hat h}_{jk})$\vspace*{1pt} is the $p\times p$ matrix obtained
when ${\hat h}_{jk}$ is defined as at (\ref{52}), and we put $\De
=\hH-H$, then, provided that $\|\De\|/\la\leq\rho$ where $\rho\in
(0,1)$ is fixed, we have
%
%e7.18 #&#
%
\begin{equation}\label{617}\qquad
\hH\mo=(I+H\mo\De)\mo H\mo
=[I-H\mo\De+O_p\{(\|\De\|/\la)^2\}] H\mo.
\end{equation}
Here the matrix $M$ represented by $O_p\{(\|\De\|/\la)^2\}$ is
interpreted as having the property $\|Mv\|\leq(1-\rho)\mo(\|\De\|
/\la)^2 \|v\|$ for all $p$-vectors $v$ (provided that $\|\De\|/\la
\leq\rho$), where on this occasion $\|Mv\|$ and $\|v\|$ denote vector
norms of the indicated quantities, and $\|\De\|$ is the matrix norm
of $\De$.

We know from (\ref{514}) that ${\hat h}_{jk}=h_{jk}+n\mhf\De
_{1jk}+n\mo\De_{2jk}$, where
%
%e7.19 #&#
%
\begin{eqnarray}\label{618}
\De_{1jk}&=&\inti\{\xi_{j+1} K^k(b)+K^{j+1}(b)
\xi_k\},\nonumber\\[-8pt]\\[-8pt]
|\De_{2jk}|&=&O_p(\th_1^j \trip K\trip^k+n\mo\trip K\trip
^{j+k}),\nonumber
\end{eqnarray}
the latter property holding uniformly in $1\leq j\leq k\leq C n\half
$. Note too that, by (\ref{610}), $\|\xi_j\|=O_p(j \th_1^j)$,
uniformly in $j\geq1$, and that $\|K^j(b)\|=O(\th_1^j)$, so $|\De
_{1jk}|=O_p\{\max(j,k) \th_1^{j+k}\}$. Therefore, if we define $\De
_{jk}={\hat h}_{jk}-h_{jk}$ then, since $\th_1<\trip K\trip<1$, we have
$n \sum\sum_{j,k\leq p} \De_{jk}^2=O_p(1)$, uniformly in $p\leq
C n\half$. Hence, $\|\De\|=O_p(n\mhf)$, uniformly in $p\leq C
n\half$, where $\De=(\De_{jk})$ is a~$p\times p$ matrix.
Therefore,
if $p$ is chosen to diverge so slowly that $p=O(n\half)$ and $\la=\la
(p)$ satisfies $n\half\la\rai$ then, by (\ref{617}),
%
%e7.20 #&#
%
\begin{equation}\label{619}
\hH\mo=\{I-H\mo\De+O_p(n\mo\la\mt)\}
H\mo,
\end{equation}
uniformly in $p\leq C n\half$. [Here $O_p(n\mo\la\mt)$ denotes a~$p\times p$ matrix, $M$ say, for which $\|Mv\|/\|v\|=O_p(n\mo\la\mt
)$ uniformly in nonzero $p$-vectors $v$.] Note too that, if we define
$\De_\ell$ to be the $p\times p$ matrix with $(j,k)$th element $\De
_{\ell jk}$, for $\ell=1,2$, then, in view of the second formula at
(\ref{618}), $\sum\sum_{j,k\leq p} \De_{2jk}^2=O_p(1)$, and
so $\|\De_2\|=O_p(1)$ uniformly\vspace*{1pt} in $p\leq C n\half$. Therefore
(\ref{619}) and the property $\De=\hH-H=n\mhf\De_1+n\mo\De
_2$ imply that
%
%e7.21 #&#
%
\begin{equation}\label{620}
\hH\mo=\{I-n\mhf H\mo\De_1+O_p(n\mo\la\mt
)\} H\mo.
\end{equation}
[Here we used the fact that $\lambda\leq h_{1,1}=O(1)$.]\vadjust{\goodbreak}%\vspace*{1pt}

Recalling the definitions of ${\hat h}_{jk}$, $\hal_j$ and $\a_j$ at
(\ref{52}), (\ref{53}) and (\ref{55}), we deduce that $\hal _j={\hat
h}_{0j}$. Noting that result (\ref{514}) can be extended to ${\hat
h}_{0j}$, we have that $\hal_j=\a_j+n\mhf\De_{10j}+n\mo \De_{20j}$,
where $\De_{10j}$ and $\De_{20j}$ are\vspace*{1pt} given by~(\ref
{618}). Note too that, by (\ref{51}) and (\ref{56}), $\hga
_j=(\hH\mo\hal)_j$ and $\ga_j=(H\mo\alpha)_j$, where $\a=(\a~_1,\ldots,\a_p)\T$ and $\hal=(\hal_1,\ldots,\hal_p)\T$. Since
$K^j(b)=O(\th_1^j)$ uniformly in $j\geq1$, $\|\eta_j\|=O_p(\trip
K\trip^j)$ uniformly in $1\leq j\leq C n\half$ (see Theorem~\ref{Theorem4}) and $\|\xi_j\|=O_p(j \th_1^j)$ uniformly in $j\geq1$ [see
(\ref{610})], then, by (\ref{510}), $\|\hkjb\|=O_p(\th _1^j+n\mhf j
\th_1^j+n\mo\trip K\trip^j)$ uniformly in $1\leq j\leq C n\half$. Using
formula (\ref{53}) for $\hal_j$, and the fact that $0<\th_1<\trip
K\trip<1$, we deduce that
%
%e7.22 #&#
%
\begin{equation}\label{623}
\|\hal\|\leq\Biggl\{\sumjop\|\hkb\|^2 \|\hkjb
\|^2\Biggr\}^{ 1/2}
=O_p(1) ,
\end{equation}
uniformly in $1\leq p\leq C n\half$.

Therefore, defining $\delta=(\De_{101},\ldots,\De_{10p})\T$, we
have, by (\ref{620}),
%
%e7.23 #&#
%
\begin{eqnarray}\label{621}
\hga&=&\hH\mo\hal\nonumber\\
&=&H\mo(\a+n\mhf\delta)
-n\mhf H\mo\De_1 H\mo\a+O_p(n\mo\la\mth
)\nonumber\\
&=&\ga+n\mhf H\mo(\delta-\De_1 \ga)
+O_p(n\mo\la\mth) ,
\end{eqnarray}
uniformly in $1\leq p\leq C n\half$, where the two vectors denoted by
$O_p(n\mo\la\mth)$ have the property that their norms equal
$O_p(n\mo\la\mth)$ uniformly in $1\leq p\leq C n\half$.

Next we combine (\ref{616}) and (\ref{621}), obtaining
%
%e7.24 #&#
%
\begin{eqnarray}\label{622}
&&\hg_p(x)-g_p(x)-(\bY-EY)\nonumber\\
&&\qquad=\sumjop\biggl[n\mhf\{H\mo(\delta-\De_1 \ga)\}
_j \inti(x-EX) K^j(b)\nonumber\\
&&\hspace*{18.5pt}\qquad\quad{}
+\ga_j\inti\{n\mhf(x-EX) \xi_j-(\bX-EX) K^j(b)\}\biggr]\nonumber\\
&&\qquad\quad{} +O_p\Biggl[n\mo\la\mth+n\mhf\sumjop
\{|\hga_j-\ga_j| (j \th_1^j+n\mhf\trip K\trip^j
)\nonumber\\
&&\qquad\quad\hspace*{138.5pt}{}
+n\mhf(|\hga_j|+|\ga_j|) \trip K\trip^j\}\Biggr] ,
\end{eqnarray}
uniformly in $1\leq p\leq C n\half$ and $\|x\|\leq C$ for each $C>0$.
Here we have used the fact that, if $V=(V_1,\ldots,V_p)\T$ is the
vector denoted by $O_p(n\mo\la\mth)$ on the far right-hand side of
(\ref{621}), then
\begin{eqnarray*}
\sumjop\biggl|V_j \inti(x-EX) K^j(b)\biggr|
&\leq&\|V\| \Biggl\{\sumjop\biggl|\inti(x-EX) K^j(b)\biggr|^2
\Biggr\}^{ 2}\\
&=&O_p(n\mo\la\mth) ,
\end{eqnarray*}
uniformly in $1\leq p\leq C n\half$ and $\|x\|\leq C$, since $\sum
_{j\geq1}\|K^j(b)\|^2<\infty$.

Note too that, since $\|K^j(b)\|=O(\th_1^j)$ and $\|\xi_j\|=O_p(j
\th_1^j)$, uniformly in $1\leq j\leq C n\half$, then by (\ref
{618}), $|\De_{1jk}|=O_p(jk \th_1^{j+k})$, uniformly in $1\leq
j,k\leq C n\half$, and therefore,
\begin{eqnarray*}
\|\De_1\|^2&\leq&\sumjop\sumkop\De_{1jk}^2=O_p(1) ,\\
\|\delta\|^2&=&\sumjop\De_{10j}^2=O_p(1) ,
\end{eqnarray*}
uniformly in $1\leq p\leq C n\half$. Hence, by (\ref{621}) and
(\ref{623}),
\begin{eqnarray*}
\|\hga-\ga\|
&=&O_p\{n\mhf\la\mo(\|\delta\|+\|\De_1\| \|\ga\|)
+n\mo\la\mth\}\\
&=&O_p\{n\mhf\la\mo(1+\|\ga\|)
+n\mo\la\mth\} .
\end{eqnarray*}
Therefore,
%
%e7.25 #&#
%
\begin{eqnarray} \label{624}
&&\sumjop
\{|\hga_j-\ga_j| (j \th_1^j+n\mhf\trip K\trip^j)
+n\mhf(|\hga_j|+|\ga_j|) \trip K\trip^j\}\nonumber\\
&&\qquad=O_p(\|\hga-\ga\|+n\mhf\|\ga\|)\nonumber\\
&&\qquad=O_p\{n\mhf\la\mo(1+\|\ga\|)
+n\mo\la\mth\} .
\end{eqnarray}
Result (\ref{515}) is a~consequence of (\ref{621}), and (\ref
{516}) follows from (\ref{622}) and~(\ref{624}).

%s7.6 #&#
\subsection{\texorpdfstring{Proof of (\protect\ref{518})}{Proof of (5.11)}}\label{proof518}
To derive (\ref{518}), note that minor modifications of the argument
used to derive (\ref{516}) can be employed to show
that, under the conditions of Theorem~\ref{Theorem5},
%
%e7.26 #&#
%
\begin{eqnarray}\label{519}
&&\|\hg_p(X_0)-g_p(X_0)\|\pred\nonumber\\
&&\qquad=\Biggl\|\bY-EY\nonumber\\
&&\qquad\quad\hspace*{3.4pt}{}+n\mhf\sumjop\biggl[\{H\mo(\delta-\De_1 \ga)\}_j
\inti(X_0-EX) K^j(b)\nonumber\\
&&\qquad\quad\hspace*{62pt}{}+\ga_j\inti\{(X_0-EX) \xi_j-n\half(\bX-EX) K^j(b)\}
\biggr]\Biggr\|\pred\nonumber\\
&&\qquad\quad{}+O_p(n\mo\la\mo\|\ga\|
+n^{-1} \la\mth) ,
\end{eqnarray}
uniformly in $p$ satisfying $1\leq p\leq C n\half$, for each $C>0$.
The predictive norm on the right-hand side of (\ref{519}) can be
shown to equal $O_p\{n\mhf\la\mo(1+\|\ga\|)\}$, and so if
(\ref{517}) holds, then
%
%e7.27 #&#
%
\begin{equation}\label{520}
\|\hg_p(X_0)-g_p(X_0)\|\pred
=O_p\{n\mhf\la\mo(1+\|\ga\|)+n\mo\la\mth\}.
\end{equation}
Since $\|g_p(X_0)-g(X_0)\|\pred=t_p(\ga_1,\ldots,\ga_p)\half$ then
(\ref{520}) implies (\ref{518}).

%apA #&#
%
\begin{appendix}
\section*{Appendix}
%sA.1 #&#
\subsection{\texorpdfstring{Conditions (\protect\ref{522}) and (\protect\ref{523})}{Conditions (5.13) and (5.14)}}\label{secA1}

Here we give examples where (\ref{522}) and (\ref{523})
hold. Assume that $E(X)=0$. Then the Karhunen--Lo\`eve expansion of
$X_i$, founded on the principal component basis introduced in
Section~\ref{sec23}, is given by $X_i=\sumj\th_j\half\xi_{ij}
\phi_j$, where the random variables $\xi_{ij}$, for $j\geq1$, are
uncorrelated and have zero mean and unit variance. For simplicity we
suppose that they have identical distributions with bounded fourth
moments, that $E(\ep^4)<\infty$, and that the eigenvalues $\th_j$
and eigenvectors $\phi_j$ satisfy the condition
$
\sumjoi\th_j\half\sup_{t\in\cI} |\phi_j(t)|<\infty.
$
Then,
%
%eA.1 #&#
%
\setcounter{equation}{0}
\begin{eqnarray}\label{A1}
&&E\Biggl[\sup_{t\in\cI} \Biggl|{1\over n\half}
\sumion\{X_i(t) D_i-EX_i(t) D_i\}\Biggr|\Biggr]\nonumber\\
&&\qquad\leq\sumjoi\th_j\half\Bigl\{\sup_{t\in\cI} |\phi_j(t)|
\Bigr\}
E\Biggl|{1\over n\half} \sumion(1-E) \xi_{ij} D_i\Biggr|
\nonumber\\
&&\qquad\leq(E\xi_{11}^4\cdot ED_1^4)^{1/4}
\sumjoi\th_j\half\Bigl\{\sup_{t\in\cI} |\phi_j(t)|\Bigr\}
<\infty,
\end{eqnarray}
%
%eA.2 #&#
%
\begin{eqnarray}
\label{A2}
&&E\Biggl[\sup_{t\in\cI} \Biggl|{1\over n\half} \sumion(1-E) \{
X_i(s)-EX_i(s)\}
\{X_i(t)-EX_i(t)\}\Biggr|^2\Biggr]\nonumber\\
&&\qquad=E\Biggl[\sup_{t\in\cI} \Biggl|\sumjoi\sumkoi(\th_j \th
_k)\half\phi_j(s) \phi_k(t)
\Biggl\{n\mhf\sumion(1-E) \xi_{ij} \xi_{ik}\Biggr\}\Biggr|^2
\Biggr]\nonumber\\
&&\qquad\leq E(\xi_{11}^4)
\Biggl[\sumjoi\th_j\half\Bigl\{\sup_{t\in\cI} |\phi_j(t)|
\Bigr\}\Biggr]^4 ,
\end{eqnarray}
where we have used the properties
\begin{eqnarray*}
\Biggl\{E\Biggl|{1\over n\half} \sumion(1-E) \xi_{ij} D_i
\Biggr|\Biggr\}^2
&\leq& E\{(\xi_{11} D_1)^2\}
\leq(E\xi_{11}^4\cdot ED_1^4)\half,\\
\Biggl\{E\Biggl|{1\over n\half} \sumion(1-E) \xi_{ij} \xi_{ik}
\Biggr|\Biggr\}^2
&\leq& E\{(\xi_{1j} \xi_{1k})^2\}
\leq E(\xi_{11}^4) .
\end{eqnarray*}
Properties (\ref{522}) and (\ref{523}) follow from (\ref
{A1}) and (\ref{A2}), respectively.

%sA.2 #&#
\subsection{Conventional implementation via the PLS basis}\label{sec31}
Inference is based on a~dataset $\cX=\{(X_1,Y_1),\ldots,
(X_n,Y_n)\}$ of independent data pairs distributed as $(X,Y)$. We first
introduce the centred data $X_i^{[1]}=X_i-\bX$ and $Y_i^{[1]}=Y_i-\bY
$, for $1\leq i\leq n$. Here and below, a~superscript in square
brackets denotes the number, or index, of a~step in our algorithm. The
algorithm goes as follows. For $j=1,\ldots,p$:

(1) Estimate $\psi_j$ by the empirical covariance of $ X_i^{[j]}$ and
$Y_i^{[j]}$:
$
\hpsi_j=\break\sumion X_i^{[j]}\times Y_i^{[j]}/\|\sumion
X_i^{[j]} Y_i^{[j]}\|.
$

(2) Fit the models
$Y_i^{[j]}=\beta_j \inti X_i^{[j]} \hpsi_j +\varepsilon_i^{[j]}$
and
$X_i^{[j]}(t)=\delta_j(t) \inti X_i^{[j]} \hpsi_j +\eta_i^{[j]}(t)$
by least-squares; that is, take
\begin{eqnarray*}
\hat\beta_j&=&\sumion Y_i^{[j]}\inti X_i^{[j]} \hpsi_j\bigg/
\sumion\biggl\{\inti X_i^{[j]} \hpsi_j\biggr\}^2 , \\
\hat\delta_j(t)&=&\sumion X_i^{[j]}(t)\inti X_i^{[j]} \hpsi_j
\bigg/\sumion\biggl\{\inti X_i^{[j]} \hpsi_j\biggr\}^2.
\end{eqnarray*}

(3)
Calculate
$X_i^{[j+1]}(t)=X_i^{[j]}(t)-\hat\delta_j(t)\inti X_i^{[j]} \hpsi
_j$ and
$Y_i^{[j+1]}=Y_i^{[j]}-\break \hat\beta_j\inti X_i^{[j]} \hpsi_j$.

After\vspace*{1pt} completion of steps (1) to (3) for all $j$, define
$M=(M_{j,k})_{1\leq j,k\leq p}$ by $M\mo=(\inti{\hat
\delta_j}\hpsi_k)_{1\leq j,k\leq p}$.
Then
$\hat b_p(t)=\sum_{j,k=1}^p \hat\beta_k M_{j,k}\hpsi_j(t)$
and
$\tg_p(x)=\bar Y+\inti\hat b_p (x-\bX) .%\label{33}
$

%sA.3 #&#
\subsection{Modified Gram--Schmidt algorithm}\label{secgram}
This algorithm turns a~set of linearly independent functions
$v_1,\ldots,v_p$ into a~set of orthogonal functions $u_1,\ldots,
u_p$, where orthogonality is defined with respect to a~scalar product
$\langle\cdot, \cdot\rangle$. For example, for the second algorithm
in Section~\ref{secaltAPLS}, the\vspace*{1pt} scalar product between two functions
$f_1$ and $f_2$ is defined by $\langle f_1 , f_2\rangle=\inti\inti
f_1(s)f_2(t) \hat K(s,\break t) \,ds \,dt$. The modified Gram--Schmidt algorithm
is described in \citet{Lan99}, Section~7.7. It works as follows:
\begin{eqnarray*}
&&\mbox{for $j=1,\ldots,p$}\\
&&\hspace*{1cm}u_j^{[1]}=v_j\\
&&\hspace*{1cm}\mbox{for $i=1,\ldots,j-1$}\\
    &&\hspace*{2cm}u_j^{[i+1]}=u_j^{[i]}-\bigl\langle u_j^{[i]},u_i\bigr\rangle  u_i\\
&&\hspace*{1cm}\mbox{end loop $i$}\\
&&\hspace*{1cm}u_j=u_j^{[j]}/\bigl\| u_j^{[j]}\bigr\|\\
&&\mbox{end loop $j$}.
\end{eqnarray*}
\end{appendix}

\section*{Acknowledgments}

We are grateful to Peter Forrester and Alan McIntosh for helpful
discussion.

%suskaldyti doi

% imsref loaded by lrinkeviciute, 2012-02-09 08:07:30
%

\printaddresses


\begin{thebibliography}{35}
% BibTex style file: ims.bst, 2011-05-30
% Default style options (sort=0,type=number).
% Used options (sort=1,type=nameyear).

%b1 #&#
\bibitem[\protect\citeauthoryear{Aguilera et~al.}{2010}]{Aguetal10}
%
\begin{barticle}[auto:STB|2012/02/03|11:55:16]
\bauthor{\bsnm{Aguilera},~\bfnm{M.}\binits{M.}},
\bauthor{\bsnm{Escabiasa},~\bfnm{M.}\binits{M.}},
\bauthor{\bsnm{Preda},~\bfnm{C.}\binits{C.}} \AND
\bauthor{\bsnm{Saporta},~\bfnm{G.}\binits{G.}}
(\byear{2010}).
\btitle{Using basis expansions for estimating functional PLS regression:
Applications with chemometric data}.
\bjournal{Chemom. Intell. Lab.}
\bvolume{104}
\bpages{289--305}.
\bptok{imsref}%
\end{barticle}
%
\endbibitem

%b2 #&#
\bibitem[\protect\citeauthoryear{Apanasovich and Goldstein}{2008}]{ApaGol08}
%
\begin{barticle}[mr]
\bauthor{\bsnm{Apanasovich},~\bfnm{Tatiyana~V.}\binits{T.~V.}} \AND
\bauthor{\bsnm{Goldstein},~\bfnm{Edward}\binits{E.}}
(\byear{2008}).
\btitle{On prediction error in functional linear regression}.
\bjournal{Statist. Probab. Lett.}
\bvolume{78}
\bpages{1807--1810}.
\bid{doi={10.1016/j.spl.2008.01.035}, issn={0167-7152}, mr={2528552}}
\bptok{imsref}%
\end{barticle}
%
\endbibitem

%b3 #&#
\bibitem[\protect\citeauthoryear{Baillo}{2009}]{Bal09}
%
\begin{barticle}[mr]
\bauthor{\bsnm{Baillo},~\bfnm{Amparo}\binits{A.}}
(\byear{2009}).
\btitle{A note on functional linear regression}.
\bjournal{J. Stat. Comput. Simul.}
\bvolume{79}
\bpages{657--669}.
\bid{doi={10.1080/00949650701836765}, issn={0094-9655}, mr={2523018}}
\bptok{imsref}%
\end{barticle}
%
\endbibitem

%b4 #&#
\bibitem[\protect\citeauthoryear{Berg and Szwarc}{2011}]{BerSzw11}
%
\begin{barticle}[mr]
\bauthor{\bsnm{Berg},~\bfnm{Christian}\binits{C.}} \AND
\bauthor{\bsnm{Szwarc},~\bfnm{Ryszard}\binits{R.}}
(\byear{2011}).
\btitle{The smallest eigenvalue of {H}ankel matrices}.
\bjournal{Constr. Approx.}
\bvolume{34}
\bpages{107--133}.
\bid{doi={10.1007/s00365-010-9109-4}, issn={0176-4276}, mr={2796093}}
\bptnote{check year}%
\bptok{imsref}%
\end{barticle}
%
\endbibitem

%b5 #&#
\bibitem[\protect\citeauthoryear{Bro and Eld{\'e}n}{2009}]{BroEld09}
%
\begin{barticle}[auto:STB|2012/02/03|11:55:16]
\bauthor{\bsnm{Bro},~\bfnm{R.}\binits{R.}} \AND
\bauthor{\bsnm{Eld{\'e}n},~\bfnm{L.}\binits{L.}}
(\byear{2009}).
\btitle{PLS works}.
\bjournal{J. Chemom.}
\bvolume{23}
\bpages{69--71}.
\bptok{imsref}%
\end{barticle}
%
\endbibitem

%b6 #&#
\bibitem[\protect\citeauthoryear{Cai and Hall}{2006}]{CaiHal06}
%
\begin{barticle}[mr]
\bauthor{\bsnm{Cai},~\bfnm{T.~Tony}\binits{T.~T.}} \AND
\bauthor{\bsnm{Hall},~\bfnm{Peter}\binits{P.}}
(\byear{2006}).
\btitle{Prediction in functional linear regression}.
\bjournal{Ann. Statist.}
\bvolume{34}
\bpages{2159--2179}.
\bid{doi={10.1214/009053606000000830}, issn={0090-5364}, mr={2291496}}
\bptok{imsref}%
\end{barticle}
%
\endbibitem

%b7 #&#
\bibitem[\protect\citeauthoryear{Cardot and Sarda}{2008}]{CarSar08}
%
\begin{barticle}[mr]
\bauthor{\bsnm{Cardot},~\bfnm{Herv{\'e}}\binits{H.}} \AND
\bauthor{\bsnm{Sarda},~\bfnm{Pascal}\binits{P.}}
(\byear{2008}).
\btitle{Varying-coefficient functional linear regression models}.
\bjournal{Comm. Statist. Theory Methods}
\bvolume{37}
\bpages{3186--3203}.
\bid{doi={10.1080/03610920802105176}, issn={0361-0926}, mr={2467760}}
\bptok{imsref}%
\end{barticle}
%
\endbibitem

%b8 #&#
\bibitem[\protect\citeauthoryear{Delaigle and Hall}{2012}]{DelHal12}
%
\begin{barticle}[auto:STB|2012/02/03|11:55:16]
\bauthor{\bsnm{Delaigle},~\bfnm{A.}\binits{A.}} \AND
\bauthor{\bsnm{Hall},~\bfnm{P.}\binits{P.}}
(\byear{2012}).
\btitle{Achieving near-perfect classification for functional data}.
\bjournal{J.~Roy. Statist. Soc. Ser. B}
\bvolume{74}
\bpages{267--286}.
\bptok{imsref}%
\end{barticle}
%
\endbibitem

%b9 #&#
\bibitem[\protect\citeauthoryear{Durand and Sabatier}{1997}]{DurSab97}
%
\begin{barticle}[mr]
\bauthor{\bsnm{Durand},~\bfnm{Jean-Fran{\c{c}}ois}\binits{J.-F.}} \AND
\bauthor{\bsnm{Sabatier},~\bfnm{Robert}\binits{R.}}
(\byear{1997}).
\btitle{Additive splines for partial least squares regression}.
\bjournal{J.~Amer. Statist. Assoc.}
\bvolume{92}
\bpages{1546--1554}.
\bid{doi={10.2307/2965425}, issn={0162-1459}, mr={1615264}}
\bptok{imsref}%
\end{barticle}
%
\endbibitem

%b10 #&#
\bibitem[\protect\citeauthoryear{Escabias, Aguilera and
Valderrama}{2007}]{EscAguVal07}
%
\begin{barticle}[mr]
\bauthor{\bsnm{Escabias},~\bfnm{M.}\binits{M.}},
\bauthor{\bsnm{Aguilera},~\bfnm{A.~M.}\binits{A.~M.}} \AND
\bauthor{\bsnm{Valderrama},~\bfnm{M.~J.}\binits{M.~J.}}
(\byear{2007}).
\btitle{Functional {PLS} logit regression model}.
\bjournal{Comput. Statist. Data Anal.}
\bvolume{51}
\bpages{4891--4902}.
\bid{doi={10.1016/j.csda.2006.08.011}, issn={0167-9473}, mr={2364547}}
\bptok{imsref}%
\end{barticle}
%
\endbibitem

\bibitem[\protect\citeauthoryear{Ferraty and Vieu}{2006}]{FerVie06}
\begin{bbook}[mr]
\bauthor{\bsnm{Ferraty},~\bfnm{Fr{\'e}d{\'e}ric}\binits{F.}} \AND
  \bauthor{\bsnm{Vieu},~\bfnm{Philippe}\binits{P.}}
(\byear{2006}).
\btitle{Nonparametric Functional Data Analysis}.
\bpublisher{Springer}, \baddress{New York}.
\bptok{imsref}%
\end{bbook}
\endbibitem

%b11 #&#
\bibitem[\protect\citeauthoryear{Frank and Friedman}{1993}]{FraFri93}
%
\begin{barticle}[auto:STB|2012/02/03|11:55:16]
\bauthor{\bsnm{Frank},~\bfnm{I.~E.}\binits{I.~E.}} \AND
\bauthor{\bsnm{Friedman},~\bfnm{J.~H.}\binits{J.~H.}}
(\byear{1993}).
\btitle{A statistical view of some chemometrics regression tools (with
discussion)}.
\bjournal{Technometrics}
\bvolume{35}
\bpages{109--148}.
\bptok{imsref}%
\end{barticle}
%
\endbibitem

%b12 #&#
\bibitem[\protect\citeauthoryear{Garthwaite}{1994}]{Gar94}
%
\begin{barticle}[mr]
\bauthor{\bsnm{Garthwaite},~\bfnm{Paul~H.}\binits{P.~H.}}
(\byear{1994}).
\btitle{An interpretation of partial least squares}.
\bjournal{J. Amer. Statist. Assoc.}
\bvolume{89}
\bpages{122--127}.
\bid{issn={0162-1459}, mr={1266290}}
\bptok{imsref}%
\end{barticle}
%
\endbibitem

%b13 #&#
\bibitem[\protect\citeauthoryear{Goutis and Fearn}{1996}]{GouFea96}
%
\begin{barticle}[mr]
\bauthor{\bsnm{Goutis},~\bfnm{Constantinos}\binits{C.}} \AND
\bauthor{\bsnm{Fearn},~\bfnm{Tom}\binits{T.}}
(\byear{1996}).
\btitle{Partial least squares regression on smooth factors}.
\bjournal{J.~Amer. Statist. Assoc.}
\bvolume{91}
\bpages{627--632}.
\bid{doi={10.2307/2291658}, issn={0162-1459}, mr={1395730}}
\bptok{imsref}%
\end{barticle}
%
\endbibitem

%b14 #&#
\bibitem[\protect\citeauthoryear{Hastie, Tibshirani and
Friedman}{2009}]{HasTibFri09}
%
\begin{bbook}[mr]
\bauthor{\bsnm{Hastie},~\bfnm{Trevor}\binits{T.}},
\bauthor{\bsnm{Tibshirani},~\bfnm{Robert}\binits{R.}} \AND
\bauthor{\bsnm{Friedman},~\bfnm{Jerome}\binits{J.}}
(\byear{2009}).
\btitle{The Elements of Statistical Learning:
Data Mining, Inference, and Prediction},
\bedition{2nd} ed.
\bpublisher{Springer}, \baddress{New York}.
\bid{doi={10.1007/978-0-387-84858-7}, mr={2722294}}
\bptok{imsref}%
\end{bbook}
%
\endbibitem

%b15 #&#
\bibitem[\protect\citeauthoryear{Helland}{1990}]{Hel90}
%
\begin{barticle}[mr]
\bauthor{\bsnm{Helland},~\bfnm{Inge~S.}\binits{I.~S.}}
(\byear{1990}).
\btitle{Partial least squares regression and statistical models}.
\bjournal{Scand. J. Stat.}
\bvolume{17}
\bpages{97--114}.
\bid{issn={0303-6898}, mr={1085924}}
\bptok{imsref}%
\end{barticle}
%
\endbibitem

%b16 #&#
\bibitem[\protect\citeauthoryear{H{\"o}skuldsson}{1988}]{Hos88}
%
\begin{barticle}[auto:STB|2012/02/03|11:55:16]
\bauthor{\bsnm{H{\"o}skuldsson},~\bfnm{A.}\binits{A.}}
(\byear{1988}).
\btitle{PLS regression methods}.
\bjournal{J. Chemom.}
\bvolume{2}
\bpages{211--228}.
\bptok{imsref}%
\end{barticle}
%
\endbibitem

%b17 #&#
\bibitem[\protect\citeauthoryear{Hou, Lascoux and Mu}{2005}]{HouLasMu05}
%
\begin{barticle}[mr]
\bauthor{\bsnm{Hou},~\bfnm{Qing-Hu}\binits{Q.-H.}},
\bauthor{\bsnm{Lascoux},~\bfnm{Alain}\binits{A.}} \AND
\bauthor{\bsnm{Mu},~\bfnm{Yan-Ping}\binits{Y.-P.}}
(\byear{2005}).
\btitle{Evaluation of some {H}ankel determinants}.
\bjournal{Adv. in Appl. Math.}
\bvolume{34}
\bpages{845--852}.
\bid{doi={10.1016/j.aam.2004.09.005}, issn={0196-8858}, mr={2129000}}
\bptnote{check year}%
\bptok{imsref}%
\end{barticle}
%
\endbibitem

%b18 #&#
\bibitem[\protect\citeauthoryear{Kr{\"a}mer, Boulesteix and
Tutz}{2008}]{KraBouTut08}
%
\begin{barticle}[auto:STB|2012/02/03|11:55:16]
\bauthor{\bsnm{Kr{\"a}mer},~\bfnm{N.}\binits{N.}},
\bauthor{\bsnm{Boulesteix},~\bfnm{A.~L.}\binits{A.~L.}} \AND
\bauthor{\bsnm{Tutz},~\bfnm{G.}\binits{G.}}
(\byear{2008}).
\btitle{Penalized partial least squares with applications to B-spline
transformations and functional data}.
\bjournal{Chemom. Intell. Lab.}
\bvolume{94}
\bpages{60--69}.
\bptok{imsref}%
\end{barticle}
%
\endbibitem

%b19 #&#
\bibitem[\protect\citeauthoryear{Kr{\"a}mer and Sugiyama}{2011}]{KraSug11}
%
\begin{barticle}[auto:STB|2012/02/03|11:55:16]
\bauthor{\bsnm{Kr{\"a}mer},~\bfnm{N.}\binits{N.}} \AND
\bauthor{\bsnm{Sugiyama},~\bfnm{M.}\binits{M.}}
(\byear{2011}).
\btitle{The degrees of freedom of partial least squares regression}.
\bjournal{J. Amer. Statist. Assoc.}
\bvolume{106}
\bpages{697--705}.
\bptok{imsref}%
\end{barticle}
%
\endbibitem

%b20 #&#
\bibitem[\protect\citeauthoryear{Lange}{1999}]{Lan99}
%
\begin{bbook}[mr]
\bauthor{\bsnm{Lange},~\bfnm{Kenneth}\binits{K.}}
(\byear{1999}).
\btitle{Numerical Analysis for Statisticians}.
\bpublisher{Springer}, \baddress{New York}.
\bid{mr={1681963}}
\bptok{imsref}%
\end{bbook}
%
\endbibitem

%b21 #&#
\bibitem[\protect\citeauthoryear{Lascoux}{1990}]{Las90}
%
\begin{barticle}[mr]
\bauthor{\bsnm{Lascoux},~\bfnm{Alain}\binits{A.}}
(\byear{1990}).
\btitle{Inversion des matrices de {H}ankel}.
\bjournal{Linear Algebra Appl.}
\bvolume{129}
\bpages{77--102}.
\bid{doi={10.1016/0024-3795(90)90299-R}, issn={0024-3795}, mr={1053054}}
\bptok{imsref}%
\end{barticle}
%
\endbibitem

%b22 #&#
\bibitem[\protect\citeauthoryear{Lorber, Wangen and
Kowalski}{1987}]{LorWanKow87}
%
\begin{barticle}[auto:STB|2012/02/03|11:55:16]
\bauthor{\bsnm{Lorber},~\bfnm{A.}\binits{A.}},
\bauthor{\bsnm{Wangen},~\bfnm{L.~E.}\binits{L.~E.}} \AND
\bauthor{\bsnm{Kowalski},~\bfnm{B.~R.}\binits{B.~R.}}
(\byear{1987}).
\btitle{A theoretical foundation for the PLS algorithm}.
\bjournal{J. Chemom.}
\bvolume{1}
\bpages{19--31}.
\bptok{imsref}%
\end{barticle}
%
\endbibitem

%b23 #&#
\bibitem[\protect\citeauthoryear{Martens and Naes}{1989}]{MarNs89}
%
\begin{bbook}[mr]
\bauthor{\bsnm{Martens},~\bfnm{Harald}\binits{H.}} \AND
\bauthor{\bsnm{Naes},~\bfnm{Tormod}\binits{T.}}
(\byear{1989}).
\btitle{Multivariate Calibration}.
\bpublisher{Wiley}, \baddress{New York}.
\bid{mr={1029523}}
\bptok{imsref}%
\end{bbook}
%
\endbibitem

%b24 #&#
\bibitem[\protect\citeauthoryear{M{\"u}ller and Yao}{2010}]{MulYao10}
%
\begin{barticle}[mr]
\bauthor{\bsnm{M{\"u}ller},~\bfnm{Hans-Georg}\binits{H.-G.}} \AND
\bauthor{\bsnm{Yao},~\bfnm{Fang}\binits{F.}}
(\byear{2010}).
\btitle{Additive modelling of functional gradients}.
\bjournal{Biometrika}
\bvolume{97}
\bpages{791--805}.
\bid{doi={10.1093/biomet/asq056}, issn={0006-3444}, mr={2746152}}
\bptok{imsref}%
\end{barticle}
%
\endbibitem

%b25 #&#
\bibitem[\protect\citeauthoryear{Nguyen and Rocke}{2004}]{NguRoc04}
%
\begin{barticle}[mr]
\bauthor{\bsnm{Nguyen},~\bfnm{Danh~V.}\binits{D.~V.}} \AND
\bauthor{\bsnm{Rocke},~\bfnm{David~M.}\binits{D.~M.}}
(\byear{2004}).
\btitle{On partial least squares dimension reduction for microarray-based
classification: A simulation study}.
\bjournal{Comput. Statist. Data Anal.}
\bvolume{46}
\bpages{407--425}.
\bid{doi={10.1016/j.csda.2003.08.001}, issn={0167-9473}, mr={2067030}}
\bptok{imsref}%
\end{barticle}
%
\endbibitem

%b26 #&#
\bibitem[\protect\citeauthoryear{Phatak and de~Hoog}{2003}]{PhadeH03}
%
\begin{barticle}[auto:STB|2012/02/03|11:55:16]
\bauthor{\bsnm{Phatak},~\bfnm{A.}\binits{A.}} \AND\bauthor{\bparticle{de}
\bsnm{Hoog},~\bfnm{F.}\binits{F.}}
(\byear{2003}).
\btitle{Exploiting the connection between PLS, Lanczos, and conjugate
gradients: Alternative proofs of some properties of PLS}.
\bjournal{J. Chemom.}
\bvolume{16}
\bpages{361--367}.
\bptok{imsref}%
\end{barticle}
%
\endbibitem

%b27 #&#
\bibitem[\protect\citeauthoryear{Phatak, Reilly and
Penlidis}{2002}]{PhaReiPen02}
%
\begin{barticle}[mr]
\bauthor{\bsnm{Phatak},~\bfnm{A.}\binits{A.}},
\bauthor{\bsnm{Reilly},~\bfnm{P.~M.}\binits{P.~M.}} \AND
\bauthor{\bsnm{Penlidis},~\bfnm{A.}\binits{A.}}
(\byear{2002}).
\btitle{The asymptotic variance of the univariate {PLS} estimator}.
\bjournal{Linear Algebra Appl.}
\bvolume{354}
\bpages{245--253}.
\bid{doi={10.1016/S0024-3795(01)00357-3}, issn={0024-3795}, mr={1927660}}
\bptok{imsref}%
\end{barticle}
%
\endbibitem

%b28 #&#
\bibitem[\protect\citeauthoryear{Preda and Saporta}{2005a}]{PreSap05N1}
%
\begin{barticle}[mr]
\bauthor{\bsnm{Preda},~\bfnm{C.}\binits{C.}} \AND
\bauthor{\bsnm{Saporta},~\bfnm{G.}\binits{G.}}
(\byear{2005}a).
\btitle{P{LS} regression on a~stochastic process}.
\bjournal{Comput. Statist. Data Anal.}
\bvolume{48}
\bpages{149--158}.
\bid{doi={10.1016/j.csda.2003.10.003}, issn={0167-9473}, mr={2134488}}
\bptok{imsref}%
\end{barticle}
%
\endbibitem

%b29 #&#
\bibitem[\protect\citeauthoryear{Preda and Saporta}{2005b}]{PreSap05N2}
%
\begin{barticle}[mr]
\bauthor{\bsnm{Preda},~\bfnm{C.}\binits{C.}} \AND
\bauthor{\bsnm{Saporta},~\bfnm{G.}\binits{G.}}
(\byear{2005}b).
\btitle{Clusterwise {PLS} regression on a~stochastic process}.
\bjournal{Comput. Statist. Data Anal.}
\bvolume{49}
\bpages{99--108}.
\bid{doi={10.1016/j.csda.2004.05.002}, issn={0167-9473}, mr={2129167}}
\bptok{imsref}%
\end{barticle}
%
\endbibitem

%b30 #&#
\bibitem[\protect\citeauthoryear{Preda, Saporta and
L{\'e}v{\'e}der}{2007}]{PreSapLev07}
%
\begin{barticle}[mr]
\bauthor{\bsnm{Preda},~\bfnm{Cristian}\binits{C.}},
\bauthor{\bsnm{Saporta},~\bfnm{Gilbert}\binits{G.}} \AND
\bauthor{\bsnm{L{\'e}v{\'e}der},~\bfnm{Caroline}\binits{C.}}
(\byear{2007}).
\btitle{P{LS} classification of functional data}.
\bjournal{Comput. Statist.}
\bvolume{22}
\bpages{223--235}.
\bid{doi={10.1007/s00180-007-0041-4}, issn={0943-4062}, mr={2318457}}
\bptok{imsref}%
\end{barticle}
%
\endbibitem

%b31 #&#
\bibitem[\protect\citeauthoryear{Reiss and Ogden}{2007}]{ReiOgd07}
%
\begin{barticle}[mr]
\bauthor{\bsnm{Reiss},~\bfnm{Philip~T.}\binits{P.~T.}} \AND
\bauthor{\bsnm{Ogden},~\bfnm{R.~Todd}\binits{R.~T.}}
(\byear{2007}).
\btitle{Functional principal component regression and functional
partial least
squares}.
\bjournal{J. Amer. Statist. Assoc.}
\bvolume{102}
\bpages{984--996}.
\bid{doi={10.1198/016214507000000527}, issn={0162-1459}, mr={2411660}}
\bptok{imsref}%
\end{barticle}
%
\endbibitem

%b32 #&#
\bibitem[\protect\citeauthoryear{Wold}{1975}]{Wol75}
%
\begin{bincollection}[mr]
\bauthor{\bsnm{Wold},~\bfnm{Herman}\binits{H.}}
(\byear{1975}).
\btitle{Soft modelling by latent variables: The non-linear iterative partial
least squares ({NIPALS}) approach}.
In \bbooktitle{Perspectives in Probability and Statistics,
Papers in Honour of M. S. Bartlett}
(\beditor{J. Gani}, ed.).
\bpublisher{Academic Press}, \baddress{London}.
\bptok{imsref}%
\end{bincollection}
%
\endbibitem

%b33 #&#
\bibitem[\protect\citeauthoryear{Wu, Fan and M{\"u}ller}{2010}]{WuFanMul10}
%
\begin{barticle}[mr]
\bauthor{\bsnm{Wu},~\bfnm{Yichao}\binits{Y.}},
\bauthor{\bsnm{Fan},~\bfnm{Jianqing}\binits{J.}} \AND
\bauthor{\bsnm{M{\"u}ller},~\bfnm{Hans-Georg}\binits{H.-G.}}
(\byear{2010}).
\btitle{Varying-coefficient functional linear regression}.
\bjournal{Bernoulli}
\bvolume{16}
\bpages{730--758}.
\bid{doi={10.3150/09-BEJ231}, issn={1350-7265}, mr={2730646}}
\bptok{imsref}%
\end{barticle}
%
\endbibitem

%b34 #&#
\bibitem[\protect\citeauthoryear{Yao and M{\"u}ller}{2010}]{YaoMul10}
%
\begin{barticle}[mr]
\bauthor{\bsnm{Yao},~\bfnm{Fang}\binits{F.}} \AND
\bauthor{\bsnm{M{\"u}ller},~\bfnm{Hans-Georg}\binits{H.-G.}}
(\byear{2010}).
\btitle{Functional quadratic regression}.
\bjournal{Biometrika}
\bvolume{97}
\bpages{49--64}.
\bid{doi={10.1093/biomet/asp069}, issn={0006-3444}, mr={2594416}}
\bptok{imsref}%
\end{barticle}
%
\endbibitem

\end{thebibliography}
\end{document}